\documentclass[a4paper, 10pt, twoside]{amsart}

\usepackage[centertags]{amsmath}
\usepackage{amsthm}
\usepackage[DIV9, headinclude, BCOR 10mm]{typearea}
\usepackage{omega_math}
\usepackage{color}
\usepackage{tabularx}
\usepackage{url}
\usepackage{mparhack} 
\usepackage{omega_draft}
\usepackage{omega_layout}
\usepackage{calc}
\usepackage{ifpdf}
\ifpdf
  \usepackage[all]{xy}
  \SelectTips{cm}{}
  \usepackage[pdftex]{graphicx}
  \usepackage[pdftex,pagebackref,colorlinks,linkcolor=blue,citecolor=blue,urlcolor=blue,a4paper,hypertexnames=false]{hyperref}
  \hypersetup{pdfauthor=Clara L\"oh,pdftitle=Isomorphisms in l\textonesuperior-homology}
  
\else
  \usepackage[all, dvips]{xy}
  \SelectTips{cm}{}
  \usepackage{graphicx}
  \newcommand{\texorpdfstring}[2]{#1}
\fi
\usepackage{omega_thm} 
\usepackage{caption}

\title{Isomorphisms in $\ell^1$-homology}
\subjclass[2000]{Primary~55N35; Secondary~46M10,~57N99}
\author{Clara L\"oh}
\date{\today}

\begin{document}

\begin{abstract}
  Taking the $\ell^1$-completion and the topological dual of the
  singular chain complex gives rise to $\ell^1$-homology and bounded
  cohomology respectively.  In contrast to \mbox{$\ell^1$-ho}\-mology,
  major structural properties of bounded cohomology are well
  understood by the work of Gromov and Ivanov.

  Based on an observation by Matsumoto and Morita,
  we derive a mechanism linking isomorphisms on the level of homology
  of Banach chain complexes to isomorphisms on the level of cohomology
  of the dual Banach cochain complexes and vice versa. Therefore,
  certain results on bounded cohomology can be transferred to
   $\ell^1$-homology. For example, we obtain a new proof of the fact
  that $\ell^1$-homology depends only on the fundamental group and
  that $\ell^1$-homology with twisted coefficients admits a
  description in terms of projective resolutions. The latter one in
  particular fills a gap in Park's approach. 

  In the second part, we demonstrate how $\ell^1$-homology can be used
  to get a better understanding of simplicial volume of non-compact
  manifolds.
\end{abstract}

\maketitle
\thispagestyle{empty}

\section{Introduction}

\noindent
Semi-norms on singular homology contain valuable geometric information
-- the fundamental example of a topological invariant created this way
is the simplicial volume of oriented, closed, connected manifolds,
which is the $\ell^1$-semi-norm of the $\R$-fundamental class.
However, singular homology itself is not an adequate algebraic tool
for the study of the $\ell^1$-semi-norm.  Only by passing to related
theories such as bounded cohomology or $\ell^1$-homology the bigger
picture becomes visible.

In contrast to $\ell^1$-homology, major structural properties of
bounded cohomology are well understood by the work of
Gromov~\cite{vbc} and Ivanov~\cite{ivanov}. For example, bounded
cohomology depends only on the fundamental group of the space in
question~\cite[p.~40, Theorem~4.3]{vbc, ivanov}, bounded cohomology
cannot see amenable normal subgroups of the fundamental
group~\cite[p.~40, Theorem~4.3]{vbc, ivanov}, and bounded cohomology
of spaces admits a description in terms of a certain flavour of
homological algebra~\cite{ivanov}.

Matsumoto and Morita observed that $\ell^1$-homology of a space is
trivial if and only if its bounded cohomology is
trivial~\cite[Corollary~2.4]{mm}. Subsequently, they raised the
natural question whether also $\ell^1$-homology depends only on the
fundamental group.  More generally one can ask how bounded cohomology
and $\ell^1$-homology are related and whether there is some kind of
duality. In the present article, we investigate to what extent such a
duality holds.

A convenient framework for this problem is the language of normed and
Banach chain complexes, i.e., chain complexes of (complete) normed
vector spaces whose boundary operators are bounded operators. Unlike
taking algebraic duals of $\R$-chain complexes, taking topological
duals of Banach chain complexes fails to commute with homology
(Section~\ref{linkingsubsec}). However, by exploiting the power of
mapping cones, we prove in Section~\ref{translproofsubsec} the
following replacement for the universal coefficient theorem:

\begin{satz}[Translation principle]\label{linkthm}\nichi{02.10.2006, 03.10.2006}
  Let $f \colon C \longrightarrow D$ be a morphism of Banach chain
  complexes and let $f' \colon D' \longrightarrow C'$ be its dual.
  \begin{enumerate}
    \item Then the induced homomorphism~$H_*(f) \colon H_*(C)
          \longrightarrow H_*(D)$ is an isomorphism of vector spaces
          if and only if $H^*(f') \colon H^*(D') \longrightarrow
          H^*(C')$ is an isomorphism of vector spaces.
    \item Furthermore, if $H^*(f') \colon H^*(D') \longrightarrow
          H^*(C')$ is an isometric isomorphism, then also $H_*(f)
          \colon H_*(C) \longrightarrow H_*(D)$ is an isometric
          isomorphism.
  \end{enumerate}
\end{satz}

In this article, the main examples for Banach (co)chain complexes are
the $\ell^1$-chain complexes and bounded cochain complexes of spaces
and of discrete groups respectively: The $\ell^1$-chain
complex~$\lonech * X$ of a topological space~$X$ is the
$\ell^1$-completion of the singular chain complex of~$X$ with
$\R$-coefficients and $\ell^1$-homology of~$X$ is defined to be the
homology of this chain complex; dually, the bounded cochain
complex~$\bchn * X$ of~$X$ is the topological dual of~$\lonech * X$
and bounded cohomology of~$X$ is defined to be the cohomology
of~$\bchn * X$. Similarly, the $\ell^1$-chain complex~$\lonech * G$ of
a discrete group~$G$ is obtained by taking the $\ell^1$-completion of
the bar resolution, and the bounded cochain complex of~$G$ is the
topological dual of~$\lonech * G$.

Applying the translation principle to suitable chain maps in the realm
of $\ell^1$-homology enables us to transfer many results concerning
bounded cohomology to $\ell^1$-homology. In particular, this strategy
provides a uniform, lightweight approach to the following results:
\begin{koro}[Isomorphisms in $\ell^1$-homology]\hfil
\begin{enumerate}
  \item Like bounded cohomology, $\ell^1$-homology of countable,
    connected \cw-complexes depends only on the fundamental group
    and amenable normal subgroups of the fundamental group are a blind
    spot of $\ell^1$-homology (Corollary~\ref{amkerkoro}).
  \item There is a characterisation of amenability of discrete groups
    through $\ell^1$-ho\-mol\-o\-gy (Corollary~\ref{amenabilitykoro}).
  \item There is a description of $\ell^1$-homology of spaces in terms
    of homological algebra; namely, $\ell^1$-homology of connected,
    countable \cw-complexes coincides with $\ell^1$-homology of the
    fundamental group, and hence $\ell^1$-homology of such spaces can
    be computed via certain strong
    relatively projective resolutions (Corollary~\ref{l1coeffcoro}).
\end{enumerate}
\end{koro}

Bouarich gave the first proof that $\ell^1$-homology depends only on
the fundamental group~\cite[Corollaire~6]{bouarich}. His proof is
based on the observation by Matsumoto and Morita, the fact that
bounded cohomology of simply connected spaces vanishes, and an
\mbox{$\ell^1$-ver}\-sion of Brown's theorem. Moreover,
Park~\cite[Corollary~4.2]{park} already claimed that
Corollary~\ref{amkerkoro} holds. However, due to a gap in her
argument, her proof is not complete. This issue is addressed in
Caveats~\ref{amenablecaveat} and~\ref{strongcaveat}, which also show
that it is not possible to imitate Ivanov's arguments in bounded
cohomology in the setting of $\ell^1$-homology.

The results listed above might give the impression that
$\ell^1$-homology is merely a shadow of bounded cohomology. However,
there are also genuine applications of \mbox{$\ell^1$-ho}\-mology: For
example, the simplicial volume of non-compact manifolds is not finite
in general -- it can even then be infinite if the manifold in question
is the interior of a compact manifold with boundary. In this case,
$\ell^1$-homology gives rise to a necessary and sufficient finiteness
condition (Theorem~\ref{fincritthm}), which cannot be phrased in terms
of bounded cohomology.

\subsection*{Organisation of this article}
In Section~\ref{normedchsec}, we introduce normed and Banach chain
complexes. 
In Section~\ref{l1bcsec}, we review the basic definitions of
$\ell^1$-homology and bounded cohomology of topological spaces as well
as of discrete groups.
Duality in the category of normed chain complexes and the proof of the
translation principle are the topic of Section~\ref{dualitysec}.
In Section~\ref{isomsec}, we apply the translation principle to
$\ell^1$-homology and we derive the consequences listed above.
Finally, in Section~\ref{fincritsec}, we demonstrate how to utilise
$\ell^1$-homology to study the simplicial volume of non-compact
manifolds. 

\subsection*{Acknowledgements}

I am grateful to Theo B\"uhler for various helpful suggestions. 

\section{Homology of normed chain complexes}\label{normedchsec}

\noindent
In this section, we introduce the basic objects of study -- normed
chain complexes and their homology. The main examples of these
concepts are $\ell^1$-homology and bounded cohomology, which are
reviewed in Section~\ref{l1bcsec}.

\subsection{Normed and Banach chain complexes}

Normed chain complexes are nothing but chain complexes in the
category of normed vector spaces (and bounded operators):

\begin{defi}[Normed chain complexes]\label{normedchdefi}\hfill 
  \begin{itemize}
    \item A \defin{normed (co)chain complex} is a (co)chain complex (indexed
      over~$\N$) consisting of normed real vector spaces, where all
      (co)boundary morphisms are bounded linear operators.
    \item A \defin{Banach (co)chain complex} is a normed (co)chain
      complex consisting of Banach spaces.
    \item A \defin{morphism of normed (co)chain complexes} is a
      (co)chain map between normed (co)chain complexes consisting of
      bounded linear operators.\qedhere
  \end{itemize} 
\end{defi}

Fundamental examples of normed chain complexes are the singular chain
complex with real coefficients and the bar resolution of a discrete
group with real coefficients (Section~\ref{l1bcsec}).

\begin{defi}[Normed chain complexes -- basic constructions] 
  Let $(C, \bou)$ be a normed chain complex.
  \begin{itemize}
    \item 
      Because the boundary operator~$\bou$ is bounded in each degree,
      it can be uniquely extended to a bounded boundary operator on the
      completion~$\overline C$ of~$C$. The resulting Banach chain
      complex, denoted by~$(\overline C, \overline \bou)$, is 
      the \defin{completion of}~$C$.
    \item 
      The \defin{dual of}~$C$ is the Banach cochain complex~$(C',
      \bou')$ defined by
      \[  (C')^n := (C_n)',
      \]
      where $\args'$ stands for the topological dual vector space,
      together with the norm given by~$\supn f := \sup \{ |f(c)| \mid c \in
      C_n,~\|c\| = 1\}$ for all~$f \in (C')^n$ and the coboundary operators 
      \begin{align*}
          (\bou')^n := (\bou_{n+1})' \colon (C')^n 
        & \longrightarrow (C')^{n+1} \\
          f
        & \longmapsto \bigl( c \mapsto f(\bou_{n+1}(c)) \bigr).
        \qedhere
      \end{align*}
  \end{itemize}
\end{defi}

Clearly, if $C$ is a normed chain complex, then $C' = (\overline
C)'$. 

\subsection{The induced semi-norm on homology}

The presence of chain complexes calls for the investigation of the
corresponding homology. In the case of normed chain complexes, the
homology groups carry additional information -- the
induced semi-norm; for example, the simplicial volume is a topological
invariant defined in terms of such a semi-norm
(Section~\ref{simvolsubsubsec}).  

\begin{defi}[Semi-norm on homology]
  Let $(C, \bou)$ be a normed chain complex, and let~$n \in \N$. The
  norm~$\|\!\cdot\!\|$ on~$C_n$ induces a semi-norm on the $n$-th homology
  group~$H_n(C) := \ker \bou_n/\im \bou_{n+1}$ as follows: If~$\alpha
  \in H_n(C)$, then
  \[     \|\alpha\| 
     :=  \inf \bigl\{ \|c\|
              \bigm|  \text{$ c\in C_n$, $\bou_n(c) = 0$, $[c] = \alpha$}
              \bigr\}.
     \qedhere
  \]
\end{defi}

In this paper, ``$\im \bou_{n+1}$'' denotes the set-theoretic image
of~$\bou_{n+1}$. 
Of course, an analogous definition applies also to normed cochain
complexes. 

Because the images of the boundary operators of a normed chain complex
are not necessarily closed, the induced semi-norm on homology in
general is not a norm; this can even happen if the underlying normed
(co)chain complex is the bounded cochain complex of a topological
space~\cite{soma,somanull}. 

Despite of the fact that the homology of a normed chain complex and
the homology of the corresponding completion in general are quite
different, the semi-norms are related. In fact, in order to understand
the semi-norms on the homology of normed chain complexes, it suffices
to consider the case of Banach chain complexes, which is shown by
approximating boundaries~\cite[Lem\-ma~2.9,
Proposition~1.7]{mschmidt,loehphd}:

\begin{prop}\label{denseprop}
  Let $D$ be a normed chain complex and let $C$ be a dense subcomplex.
  Then the map~$H_*(C) \longrightarrow H_*(D)$ induced by the
  inclusion is isometric.
\end{prop}

Moreover, one can also compute the induced semi-norm on~$H_*
(C)$ via the semi-norm on~$H^*(C')$ (Theorem~\ref{dualitygromovthm}).

%
\section{\texorpdfstring{$\ell^1$}{l\textonesuperior}-Homology and bounded
  cohomology}\label{l1bcsec}

\noindent
Taking the completion and the topological dual of the singular chain
complex with respect to the $\ell^1$-norm gives rise to
$\ell^1$-homology and bounded cohomology respectively
(Section~\ref{l1bcspacessubsec}). Also the bar resolution
of a discrete group admits an $\ell^1$-norm -- leading to
$\ell^1$-homology and bounded cohomology of discrete groups
(Section~\ref{l1bcgroupssubsec}). Both constructions can be decorated
with equivariant Banach modules, which yields the corresponding
theories with (twisted) coefficients (Sections~\ref{l1bcgroupssubsec}
and~\ref{l1bctwistedsubsec}). 

\subsection{\texorpdfstring{$\ell^1$}{l\textonesuperior}-Homology and bounded
  cohomology of spaces}
  \label{l1bcspacessubsec} 

We start with the key example of a normed chain complex:

\begin{defi}[$\ell^1$-Norm on the singular chain complex]
  Let $(X,A)$ be a pair of topological spaces. 
  \begin{itemize}
    \item 
      The \defin{$\ell^1$-norm} on the singular chain
      complex~$\csing * X$ with real coefficients is defined as
      follows: For a chain~$c = \sum_{j=0}^k a_j \cdot \sigma_j \in
      \csing n X$ in reduced form we set
      \[ \lone c := \sum_{j=0}^k |a_j|. 
      \]
    \item The induced semi-norm on the quotient~$\csing * {X,A} = \csing * X /
      \csing * A$ is a norm because the subcomplex~$\csing * A$ is
      $\ell^1$-closed in~$\csing * X$; this norm on~$\csing * {X,A}$
      is also denoted by~$\lone{\args}$.  \qedhere
  \end{itemize}
\end{defi}

The boundary operator~$\bou_n \colon \csing n {X,A} \longrightarrow
\csing {n-1} {X,A}$ is a bounded operator with respect to the
$\ell^1$-norm of operator norm at most~$(n+1)$. Hence, $\csing * {X,A}$
is a normed chain complex.  Clearly, $\csing * X$ and~$\csing * {X,A}$
are in general not complete and thus these complexes are no Banach
chain complexes.

On the other hand, for~$p \in (1,\infty]$, the singular chain complex
equipped with the $\ell^p$-norm is in general \emph{not} a normed
chain complex in the sense of
Definition~\ref{normedchdefi}~\cite[Proposition~2.11]{loehphd}.

\begin{defi}[$\ell^1$-Homology and bounded cohomology of spaces]\label{l1bcdefi}
  Let $(X,A)$ be a pair of topological spaces. 
  \begin{itemize}
    \item The \defin{$\ell^1$-chain complex of~$(X,A)$} is the 
      completion~$\lonech * {X,A}$ of the normed chain 
      complex~$\csing * {X,A}$ with respect to~$\lone\args$.
      We abbreviate~$\lonech * {X,\emptyset}$ by~$\lonech * X$. 
    \item Then \defin{$\ell^1$-homology of~$(X,A)$} is defined as
      \[    \lonehom * {X,A} 
         := H_*\bigl( \lonech * {X,A} \bigr).
      \]
    \item Dually, the \defin{bounded cochain complex of~$(X,A)$} is
      the dual~$\bchn *{X,A}$ of the normed chain complex~$\csing *
      {X,A}$. 
    \item \defin{Bounded cohomology of~$(X,A)$} is given by
      \[    \bch * {X,A} 
         := H^*\bigl( \bchn * {X,A}\bigr).
      \]
    \item The semi-norms on $\ell^1$-homology and bounded cohomology
      induced by~$\lone\args$ and~$\supn\args$ respectively are also
      denoted by~$\lone\args$ and~$\supn\args$.
      \qedhere
  \end{itemize}
\end{defi}

The inclusion~$\csing * {X,A} \hookrightarrow \lonech * {X,A}$ of
chain complexes induces a \emph{comparison map} $\hsing * {X,A}
\longrightarrow \lonehom * {X,A}$, which is isometric by
Proposition~\ref{denseprop}; in general, this homomorphism is neither
injective nor surjective. Similarly, there is a comparison map~$\bch *
{X,A} \longrightarrow \cohsing * {X,A}$.

\subsubsection{Functoriality}
If $f \colon (X,A) \longrightarrow (Y,B)$ is a continuous map of pairs
of topological spaces, then the induced map~$\csing * f \colon \csing
* {X,A} \longrightarrow \csing * {Y,B}$ is a morphism of normed chain
complexes. Consequently, we obtain induced morphisms~$\lonech * f$
and~$\bchn * f$, as well as maps~$\lonehom * f \colon \lonehom *
{X,A} \longrightarrow \lonehom * {Y,B}$ and~$\bch * f \colon
\bch * {X,A} \longrightarrow \bch * {Y,B}$; clearly, both~$\lonehom *
\args$ and~$\bch * \args$ are functorial with respect to composition.

\subsubsection{Basic properties}
Standard arguments show that $\ell^1$-homology and bounded cohomology
are homotopy invariant and admit a long exact sequence for pairs of
topological spaces~\cite[Proposition~2.7]{loehphd}. Using self-maps of
the circle of non-trivial degree one finds that $\lonehom 1 X = 0$
and~$\bch 1 X = 0$ holds for all spaces~$X$~\cite[Corollary~2.7,
Proposition~2.7]{mm, loehphd}.  However, both $\ell^1$-homology and
bounded cohomology do \emph{not} satisfy excision~\cite{brooks,
  mitsumatsu} (infinite chains need not contain only
small simplices after a finite number of barycentric subdivisions).
This failure of excision is both a curse and a blessing. On the one
hand, the lack of excision makes concrete computations via the usual
divide and conquer approach significantly harder; on the other hand,
it turns out that bounded cohomology and $\ell^1$-homology depend
only on the fundamental group and hence can be computed in terms of
certain nice resolutions (Corollary~\ref{amkerkoro} and
Corollary~\ref{l1coeffcoro}).

%
\subsubsection{Simplicial volume}\label{simvolsubsubsec}

An example of valuable geometric information encoded in a semi-norm on
homology is the simplicial volume introduced by Gromov~\cite{vbc}. The
simplicial volume is a homotopy invariant linked to Riemannian
geometry in various ways and can be viewed as a topological
approximation of the Riemannian volume~\cite{vbc}.

\begin{defi}
  Let $M$ be an oriented, closed, connected $n$-manifold with
  $\R$-fun\-da\-men\-tal class~$\fcl M \in \hsing n M$. Then the
  \defin{simplicial volume of}~$M$ is defined as
  \[    \sv M 
     := \lone[big]{\fcl M}
     =  \inf\,\bigl\{ \lone c
            \bigm|  \text{$c \in \csing n M$ is an $\R$-fundamental
                    cycle of~$M$}
            \bigr\}.
     \qedhere
  \] 
\end{defi}

Using self-maps of non-trivial degree one sees that the simplicial
volume of spheres and tori is zero. On the other hand, straightening
simplices to geodesic simplices shows that the simplicial volume of
closed hyperbolic manifolds is non-zero~\cite{thurston, inoue}. 

However, it is in general very difficult to compute the simplicial
volume by geometric means. In view of Proposition~\ref{denseprop} and
Theorem~\ref{dualitygromovthm} below and the comparison maps, it is
possible to use \mbox{$\ell^1$-ho}\-mology and bounded cohomology to
compute the simplicial volume. For example, this approach shows that
the simplicial volume of all manifolds with amenable fundamental group
is zero. Conversely, we can deduce that $\ell^1$-homology and bounded
cohomology of closed hyperbolic manifolds are non-trivial.

%
\subsection{\texorpdfstring{$\ell^1$}{l\textonesuperior}-Homology and bounded
  cohomology of discrete groups}\label{l1bcgroupssubsec}

For a discrete group~$G$, we write~$C_*(G)$ for the corresponding bar
resolution with real coefficients; more explicitly, $C_n(G)$ is the
free $\R G$-module with basis~$([g_1 | \dots | g_n])_{g\in G^{n}}$,
and the boundary operator~$C_n(G) \longrightarrow C_{n-1}(G)$ is the
$G$-linear map determined uniquely by
\begin{align*}    
  C_n(G) \longrightarrow & \; C_{n-1}(G) \\
  [g_1 | \dots | g_n ]
         \longmapsto & \;
         g_1 \cdot [g_2 | \dots | g_n]
         \\ + & \; 
         \sum_{j=1}^{n-1} (-1)^j \cdot
         [g_1 | \dots |g_{j-1}| g_j \cdot
         g_{j+1}|g_{j+2}| \dots |g_n]
         \\ + & \;
         (-1)^n \cdot [g_1 | \dots | g_{n-1}]
         .
\end{align*}

\begin{defi}[$\ell^1$-Norm on the bar resolution of discrete groups]
  Let $G$ be a discrete group, and let $n \in \N$. For $c = \sum_{g
    \in G^{n+1}} a_g \cdot g_0 \cdot [g_1 | \dots | g_n] \in C_n(G)$ we define
  \[ \lone c 
     := \sum_{g\in G^{n+1}} |a_g|.
     \qedhere
  \]
\end{defi}

The group~$G$ acts isometrically on~$C_*(G)$ and $C_*(G)$ is a normed
chain complex with respect to the $\ell^1$-norm; in particular, we obtain
the corresponding completions and topological duals:

\begin{defi}[$\ell^1$-chains and bounded cochains of discrete groups] 
  Let $G$ be a discrete group.
  \begin{itemize}
    \item The \emph{$\ell^1$-chain complex of~$G$} is the completion~$\lonech
      * G$ of the normed chain complex~$C_*(G)$ with respect
      to~$\lone{\args}$. 
    \item The \emph{bounded cochain complex of~$G$} is the dual~$\bchn
      * G$ of the normed chain complex~$C_*(G)$.  
      \qedhere
  \end{itemize}
\end{defi}

In order to define $\ell^1$-homology and bounded cohomology of
discrete groups (with coefficients), we need some terminology from the
category of Banach $G$-modules: A \emph{Banach $G$-module} is a Banach
space equipped with an isometric (left) $G$-action. If $U$ and~$V$ are
two Banach $G$-modules, then the projective tensor product~$U \projten
V$ and the space~$B(U,V)$ of bounded linear functions from~$U$ to~$V$
are Banach $G$-modules with respect to the following, diagonal,
$G$-actions: For all~$g \in G$ one sets
\begin{align*}
   \fa{u \in U} \fa{v \in V}
   g \cdot (u \otimes v) & := (g \cdot u) \otimes (g \cdot v), \text{~and}\\
   \fa{f \in B(U,V)}
   g \cdot f & := \bigl(u \mapsto g \cdot f(g^{-1} \cdot u) \bigr).
\end{align*}
For a Banach $G$-module~$V$ the set of \emph{invariants} of~$V$ is
defined by
\[ V^G := \{ v \in V
          \mid \fa{g \in G} g \cdot v = v
          \};
\]
the set of \emph{coinvariants} of~$V$ is the quotient
$ V_G := V/\overline W,
$
where $W \subset V$ is the subspace generated by the set~$\{ g \cdot v
- v \mid v \in V, g\in G \}$. It is not difficult to see that there is
an isometric isomorphism~$(V_G)' \cong (V')^G$. 

A \emph{Banach $G$-(co)chain complex} is a normed (co)chain complex
consisting of Banach $G$-modules whose (co)boundary operators are
$G$-equivariant. For example, $\lonech * G$ is a Banach $G$-chain
complex.  A \emph{morphism of Banach $G$-(co)chain complexes} is just
a morphism of normed (co)chain complexes that is $G$-equivariant.
Notions such as the invariants etc.\ have obvious analogues on the
level of Banach $G$-(co)chain complexes.

Now the definition of $\ell^1$-homology and bounded cohomology of
discrete groups is a straightforward adaption of the definition of
group (co)homology in terms of the bar resolution:

\begin{defi}[$\ell^1$-Homology and bounded cohomology of discrete
  groups]
  Let $G$ be a discrete group, and let $V$ be a Banach $G$-module.
  \begin{itemize}
    \item We write
          \[ \lonech * {G ; V} := \lonech * G \projten V 
             \quad\text{and}\quad
             \bchn * {G;V} := B\bigl(\lonech * G, V\bigr).
          \]
    \item The \emph{$\ell^1$-homology of~$G$ with coefficients
        in~$V$}, denoted by~$\lonehom * {G;V}$, is the homology of the
      Banach chain complex~$\lonech * {G;V}_G$. 
    \item \emph{Bounded cohomology of~$G$ with coefficients in~$V$},
      denoted by~$\bch * {G;V}$, is the cohomology of the Banach
      cochain complex~$\bchn * {G;V}^G$.
      \qedhere
  \end{itemize}
\end{defi}

Notice that $\bchn * {G;V'}$ is isometrically $G$-isomorphic
to~$(\lonech * {G;V})'$; in particular, we have $\bchn * {G; \R} =
\bchn * G$, where $\R$ is equipped with the trivial $G$-action. For
brevity, we write $\lonehom * G := \lonehom * {G;\R}$ and~$\bch * G :=
\bch * {G; \R}$.

Moreover, the $\ell^1$-norm on~$\lonech * G$ and the norm on~$V$
induce norms on~$\lonech * {G;V}$ and~$\bchn * {G;V}$, and hence they
induce semi-norms on~$\lonehom *{G;V}$ and~$\bch * {G;V}$. These
semi-norms are also denoted by~$\lone{\args}$ and~$\supn{\args}$
respectively.

\subsubsection{$\ell^1$-Homology and bounded cohomology in degree~$0$}\label{zerosubsubsec}

Almost the same calculations as in ordinary group (co)homology show
that $\bch 0 {G;V} \cong V^G$ and $\lonehom 0 {G;V} \cong V/U$ for all
discrete groups~$G$ and all Banach $G$-modules~$V$; here,
\[ U = 
   \biggl\{ \sum_{j \in \N} a_j \cdot (v_j - g_j \cdot v_j)
   \biggm| \text{$(a_j)_j \subset \R$, 
                 $(g_j)_j \subset G$, 
                 $(v_j)_j \subset V$
                 and $\sum_{j \in \N}|a_j| \cdot \|v_j\| < \infty$}
   \biggr\}.
\]
We have $V/\overline U = V_G$, but in general $U$ is not closed in~$V$
and so $V/U$ need not be equal to~$V_G$. If $V$ is a reflexive Banach
space, then indeed $\lonehom 0 {G;V} \cong V_G$: If $V$ is reflexive,
then $0 = \bch 1 {G;V'} \cong H^1(\lonech *
{G;V}')$~\cite[Propositon~6.2.1]{monod}. Therefore, $\lonehom 0 {G;V}$
is Banach~\cite[Theorem~2.3]{mm} and hence $\lonehom 0 {G;V} \cong
V/\overline U = V_G$.

\subsubsection{Functoriality}

Let $\varphi \colon G \longrightarrow H$ be a homomorphism of discrete
groups, let $V$ be a Banach $G$-module and let $W$ be a Banach
$H$-module.  Then 
  \begin{align*}
      \lonech n \varphi
      \colon \lonech n G 
    & \longrightarrow 
      \varphi^*\bigl(\lonech n H\bigr) \\
      g_0 \cdot [g_1 | \dots | g_n]
    & \longmapsto
      \varphi(g_0) \cdot
      \bigl[ \varphi(g_1) | \dots | \varphi(g_n)
      \bigr]
  \end{align*}
  defines a morphism~$\lonech * \varphi \colon \lonech * G
  \longrightarrow \varphi^*\lonech * H$ of Banach $G$-chain complexes
  of norm~$1$; here, $\varphi^*(\args)$ stands for the Banach
  $G$-module structure on the Banach \mbox{$H$-mod}\-ule in question that is
  induced by~$\varphi$. 
  In particular, for any morphism~$f \colon V
  \longrightarrow \varphi^* W$ of Banach $G$-modules, the map 
  \begin{align*}
    \lonech * {\varphi; f}
    := \lonech * \varphi \projten f
    \colon \lonech * {G;V}
    \longrightarrow 
    \varphi^*\bigl(\lonech * {H; W}\bigr)
  \end{align*}
  is a morphism of Banach $G$-chain complexes (of norm at
  most~$\|f\|$). Analogously, for any morphism~$f \colon \varphi^*W
  \longrightarrow V$ of Banach $G$-modules,
  \begin{align*}
    \bchn * {\varphi; f}
    := B\bigl( \lonech * \varphi, f
        \bigr)
    \colon \varphi^*\bigl(\bchn * {H; W}\bigr)
    \longrightarrow 
    \bchn * {G;V}
  \end{align*}
  is a morphism of Banach $G$-cochain complexes (of norm at
  most~$\|f\|$).

Let $p \colon ( \varphi^* \lonech * {H;W})_G \longrightarrow
\lonech * {H;W} _H$ and $i \colon \bchn * {H;W} ^H \longrightarrow
(\varphi^* \bchn * {H;W})^G$ denote the canonical projection and the
inclusion respectively. Then we write
\begin{align*}
        \lonehom * {\varphi; f} 
  & :=  H_* \bigl( p \circ \lonech * {\varphi; f}_G
            \bigr)
    \colon \lonehom * {G;V} \longrightarrow \lonehom * {H;W},
  \\
        \bch * {\varphi; f} 
  & :=  H^* \bigl( \bchn * {\varphi; f}^G \circ i
            \bigr)
    \colon \bch * {H;W}\longrightarrow \bch * {G;V}.
\end{align*}

\subsubsection{Strong relatively injective/relatively projective
  resolutions}

Both $\ell^1$-homology and bounded cohomology of discrete groups enjoy
the same flexibility as ordinary group (co)homology: namely, both
theories can be computed by means of relative homological algebra as
studied by Brooks, Ivanov, Monod, and Park~\cite{brooks, ivanov,
  monod, park}. 

As in the classical case, there is a distinguished class of
resolutions -- so-called strong relatively projective resolutions and
strong relatively injective resolutions -- and a corresponding
fundamental lemma of homological algebra granting existence and
uniqueness of certain morphisms of Banach $G$-chain
complexes~\cite[Appendix~A]{loehphd}; for example, the Banach
(co)chain complexes~$\lonech * {G;V}$ and~$\bchn * {G;V}$ together
with the obvious augmentation maps are strong
relatively projective/injective $G$-resolutions
of~$V$~\cite[Proposition~2.19]{loehphd}. Therefore, we
obtain~\cite[Theorem~2.18]{loehphd}:

\begin{satz}\label{discreteresthm}
  Let $G$ be a discrete group and let $V$ be a Banach $G$-module. 
   \begin{enumerate}
    \item 
      For any strong relatively projective $G$-resolution~$(C,
      \eta \colon C_0 \rightarrow V)$ of~$V$ there is a
      canonical isomorphism (degreewise isomorphism of semi-normed 
      vector spaces)
      \[ \lonehom * {G;V} \cong H_*(C_G).
      \]
    \item
      For any strong relatively injective $G$-resolution~$(C,
      \eta \colon V \rightarrow C^0)$ of~$V$ there is a
      canonical isomorphism (degreewise isomorphism of semi-normed
      vector spaces) 
      \[ \bch * {G;V} \cong H^*(C^G).
      \] 
    \item 
      If $(C, \eta \colon C_0 \rightarrow \R)$ is a strong
      relatively projective $G$-resolution of the trivial Banach
      $G$-module~$\R$, then there are canonical isomorphisms
      (degreewise isomorphisms of semi-normed vector spaces)
      \begin{align*}
                \lonehom * {G;V}
        & \cong H_*\bigl( (C \projten V)_G
                   \bigr), \\
                \bch * {G;V}
        & \cong H^*\bigl( B(C,V)^G
                   \bigr).
      \end{align*}
  \end{enumerate}
\end{satz}

The semi-norms on~$\lonehom * {\args;\cdot\,}$ and~$\bch *
{\args;\cdot\,}$ induced by the bar resolutions~$\lonech *
{\args;\cdot\,}$ and~$\bchn *{\args;\cdot\,}$ coincide with the
canonical semi-norms in the sense of Ivanov~\cite[Corollary~3.6.1,
Corollary~2.3, Corollary~7.4.7]{ivanov, park, monod}. On the other
hand, rescaling augmentation maps shows that not every strong relatively
projective/injective resolution induces the same semi-norm in
(co)homology.

B\"uhler developed a description of $\ell^1$-homology and bounded
cohomology as derived functors via exact categories~\cite{buehler},
thereby providing an even more conceptual approach.

\subsection{$\ell^1$-Homology and bounded cohomology of spaces with
  twisted coefficients}\label{l1bctwistedsubsec}

Similarly to singular homology and singular cohomology there are also
versions of $\ell^1$-homology and bounded cohomology of spaces with
twisted coefficients:

\begin{defi}[$\ell^1$-Homology and bounded cohomology with twisted coefficients]
  Let $X$ be a connected topological space with fundamental group~$G$
  that admits a universal covering space~$\ucov X$, and let $V$ be a
  Banach $G$-module.
  \begin{itemize}
    \item The \emph{$\ell^1$-chain complex of~$X$ with twisted
        coefficients in~$V$} is defined as the Banach chain complex of
        coinvariants
        \[ \lonech * {X; V} := \bigl( \lonech * {\ucov X} \projten V \bigr)_G. 
        \]
        Here, $\lonech* {\ucov X}$ inherits the $G$-action from the
        action of the fundamental group on the universal
        covering~$\ucov X$.
    \item The \emph{$\ell^1$-homology of~$X$ with twisted coefficients
        in~$V$}, denoted by~$\lonehom * {X;V}$, is the homology of the
        Banach chain complex~$\lonech * {X;V}$. 
    \item The \emph{bounded cochain complex of~$X$ with twisted
        coefficients in~$V$} is defined as the Banach cochain complex
        of invariants
        \[ \bchn * {X;V} := B\bigl(\lonech * {\ucov X} , V\bigr)^G. \]
    \item \emph{Bounded cohomology of~$X$ with twisted coefficients
        in~$V$} is the cohomology of the Banach cochain complex~$\bchn
        * {X;V}$ and is denoted by~$\bch * {X;V}$.
       \qedhere
  \end{itemize}
\end{defi}

The $\ell^1$-chain complex and the bounded cochain complex of~$X$ as
defined in Definition~\ref{l1bcdefi} can be recovered from this
definition by taking $\R$ with the trivial $G$-action as 
coefficients~\cite[Proposition~2.23]{loehphd}.

\section{Duality}\label{dualitysec}

\noindent
In this section, we investigate the relation induced by the evaluation
map between homology of a normed chain complex and cohomology of its
dual cochain complex. Unlike taking algebraic duals of $\R$-chain
complexes, taking topological duals of normed chain complexes fails to
commute with homology (Section~\ref{linkingsubsec}).
Section~\ref{translproofsubsec} is devoted to the proof of the
translation principle (Theorem~\ref{linkthm}), showing that it is
still possible to transfer certain information from homology of a
Banach chain complex to cohomology of the dual complex and vice versa.

\subsection{Linking homology and cohomology}\label{linkingsubsec}

Evaluation links homology of a normed chain complex to cohomology of
its dual cochain complex: If $C$ is a normed chain complex and $n \in
\N$, then the evaluation map~$C'{}^n \otimes C_n \longrightarrow \R$
induces a linear map
\[ \krp\args{\!\args} \colon H^n(C') \otimes H_n(C) \longrightarrow \R,
\]
the so-called \emph{Kronecker product}. Similarly, we obtain a
map~$\overline H{}^n(C') \longrightarrow (\overline H_n(C))'$,
where $\overline H$ denotes reduced (co)homology, i.e., the kernel
modulo the \emph{closure} of the image of the (co)boundary operator.

Taking the algebraic dual is compatible with taking homology: For all
$\R$-chain complexes~$C$ the map~$H^*\bigl(\hom_\R(C,\R)\bigr)
\longrightarrow \hom_\R \bigl( H_*(C), \R \bigr)$ induced by
evaluation is an isomorphism by the universal coefficient
theorem. However, taking topological duals, even of complete normed
chain complexes, fails to commute with taking homology:

\begin{bem}\nichi{??.10/11.2006}\label{obviousbem}
  \emph{There is no obvious duality isomorphism between homology and
  cohomology of Banach chain complexes}:

  Let $C$ be a Banach chain complex. Then we have the following
  commutative diagram 
  \begin{align*}
   \xymatrix@=2em{%
       H^*(C') \ar[r] \ar[d] \ar[dr]
     & \hom_\R( H_*(C), \R) \\
       \overline H ^*(C') \ar[r]
     & \bigl( \overline H_*(C) \bigr)', \ar[u]
   }
  \end{align*} 
  where the horizontal arrows are the homomorphisms induced by the
  Kronecker products (i.e., they are induced by evaluation of elements
  in~$C'$ on elements in~$C$), the left vertical arrow is the
  canonical projection and the right vertical arrow is the
  composition~$(\overline H_*(C))' \hookrightarrow \hom_\R(\overline
  H_*(C), \R) \hookrightarrow \hom_\R(H_*(C), \R)$ of inclusions.

  The lower horizontal morphism, and hence also the diagonal morphism,
  is surjective by the Hahn-Banach theorem. Moreover, Matsumoto and
  Morita showed that the diagonal morphism is injective if and only
  if~$H^*(C') = \overline H^*(C')$ holds~\cite[Theorem~2.3]{mm}.

  Obviously, this is not the case in general. It is even wrong if~$C =
  \lonech * X$ for certain topological spaces~$X$ \cite{soma,
  somanull}.  Hence, there is no obvious duality between
  $\ell^1$-homology and bounded cohomology.

  In addition, the lower horizontal arrow is in general not injective: The
  kernel of the evaluation map
  \[ \ker \bou'^{n+1} \longrightarrow 
     \bigl( \ker \bou_n / \overline{\im \bou_{n+1}}\bigr)'
     = \bigl( \overline H_n(C)
       \bigr)'
  \]
  equals~$({}^\bot\im (\bou'^n))^{\bot}$, which is the weak*-closure
  of~$\im \bou'^n$~\cite[Theorem~4.7]{rudin}. Furthermore, the
  norm-closure~$\overline{\im \bou'^n}$ and the
  weak*-closure~$({}^\bot\im (\bou'^n))^{\bot}$ coincide if and only
  if $\im \bou_{n+1}$ is closed~\cite[Theorem~4.14]{rudin}.  Thus
  there is also no obvious duality isomorphism between reduced
  $\ell^1$-homology and reduced bounded cohomology.
\end{bem}

Nevertheless, the Kronecker product is strong enough to give
sufficient conditions for (co)homology classes to be non-trivial. For
example, if $\alpha \in H_*(C)$ and $\varphi \in H^*(C')$
with~$\krp\varphi\alpha = 1$, then neither~$\alpha$, nor~$\varphi$ can
be zero. This effect can be used to show that $\ell^1$-homology and
bounded cohomology of certain surface groups are
non-trivial~\cite{mitsumatsu}.

\subsection{Transferring isomorphisms -- proof of Theorem~\ref{linkthm}}\label{translproofsubsec}

\subsubsection{Method of proof}

The proof of the translation principle (Theorem~\ref{linkthm}) relies
on the following three tools: 

\begin{enumerate}
  \item \emph{Duality principle.} 
        There is the following relation between homology of Banach chain
        complexes and cohomology of their duals, which has been discovered
        by Johnson as well as by Matsumoto and
        Morita~\cite[Proposition~1.2, Corollary~2.4, Theorem~3.5]{johnson,
        mm, loehphd}.

        \begin{satz}[Duality principle]\label{mmthm}
          Let $C$ be a Banach chain complex. Then $H_*(C)$ vanishes if and
          only if~$H^*(C')$ vanishes.
        \end{satz}

        Here, the ``$*$'' carries the meaning ``\emph{All} of
        the~$H_n(C)$ are zero if and only if \emph{all} of
        the~$H^n(C')$ are zero.'' 
        The key to lifting this duality principle to morphisms of
        Banach chain complexes is to apply the duality principle to
        the mapping cone of the morphism in question.
        
  \item \emph{Mapping cones.}
        Mapping cones of chain maps are a device translating questions
        about isomorphisms on homology into questions about the
        vanishing of certain homology groups; the exact definition of
        mapping cones in the context of Banach chain complexes is
        given in Section~\ref{conessubsubsec} below.

        \begin{prop}\label{conesprop}
          Let $f \colon C \longrightarrow D$ be a morphism of normed
          chain complexes. 
          \begin{enumerate}
            \item The induced map~$H_*(f) \colon H_*(C)
              \longrightarrow H_*(D)$ is an isomorphism of vector
              spaces if and only if~$H_*(\cone f) = 0$. Of course, the
              analogous statement for morphisms of normed cochain
              complexes also holds. 
            \item There is a natural isomorphism
              $\cone f ' \cong \Sigma \cone {-f'} 
              $
              of normed cochain complexes, relating the mapping cones
              of~$f$ and~$-f'$. 
          \end{enumerate}
        \end{prop}

        The suspension~$\Sigma$ just shifts the (co)chain complex in
        question by~$+1$ and changes the sign of the boundary operator.

        The first part of Proposition~\ref{conesprop} is a classic
        fact from homological algebra (long exact homology sequence
        associated with the mapping cone~\cite[Section~1.5]{weibel}); a
        straightforward calculation proves the second part.

  \item \emph{Duality principle for semi-norms.}
        The third ingredient for the proof of the translation principle is
        the following observation of Gromov~\cite[p.~17,
        Proposition~F.2.2, Theorem~3.8]{vbc,bp, loehphd}, relating the
        semi-norm on homology to the semi-norm on cohomology of the dual.

        \begin{satz}[Duality principle for semi-norms]\label{dualitygromovthm}
          Let $C$ be a normed chain complex and let~$n \in \N$. Then 
          \[ \|\alpha \| 
             = \sup
               \Bigl\{ \frac 1 {\supn \varphi}
               \Bigm|  \text{$\varphi \in H^n(C')$ and 
                             $\krp \varphi \alpha = 1$}
               \Bigr\}
          \]
          holds for each~$\alpha \in H_n(C)$; here, $\sup \emptyset := 0$.
        \end{satz}

        However, the semi-norm on cohomology of the dual can in
        general not be computed in terms of the semi-norm on homology,
        because it can happen that the reduced homology~$\overline
        H_*(C)$ is zero while $\overline H^*(C')$ is non-zero
        (cf.~Remark~\ref{obviousbem}). 

\end{enumerate}

\subsubsection{Mapping cones}\label{conessubsubsec}

        For the sake of completeness, we recall the definition of
        mapping cones of morphisms of Banach chain complexes: 

        \begin{defi}[Mapping cones]\hfill
           \begin{itemize}
           \item
             Let $f \colon (C, \bou^C) \longrightarrow (D, \bou^D)$ be
             a morphism of normed chain complexes. Then the
             \defin{mapping cone} of~$f$, denoted by~$\cone{f}$, is
             the normed chain complex defined by
             \[ \cone{f}_n := C_{n-1} \oplus D_n, \]
             linked by the boundary operator that is given by the matrix
             \begin{align*}
               \xymatrix@=3em@C=-.2em{%
                 \cone f _n     \ar_{\footnotesize\begin{pmatrix}
                                     -\bou^C         & 0 \\
                                     \phantom{-}f    & \bou^D
                                     \end{pmatrix}}[d]
                            & = C_{n-1}\; \oplus   
                                \ar_{-\bou^C}[d] \ar_{f}[dr]
                            & D_n  
                                \ar^{\bou^D}[d]\\
                              \cone f _{n-1} & = C_{n-2}\; \oplus
                            & D_{n-1}.
               }
             \end{align*}
           \item
             Dually, if $f \colon (D, \cobou_D) \longrightarrow (C,
             \cobou_C)$ is a morphism of normed cochain complexes,
             then the \defin{mapping cone} of~$f$, also denoted
             by~$\cone f$, is the normed cochain complex defined by
             \[ \cone f ^n := D^{n+1} \oplus C^{n} \]
             with the coboundary operator determined by the matrix 
             \[
               \xymatrix@=3em@C=-.2em{%
                  \cone f ^n     \ar_{\footnotesize\begin{pmatrix}
                                      -\cobou_D       & 0 \\
                                      \phantom{-}f    & \cobou_C
                                      \end{pmatrix}}[d]
                             & = D^{n+1}\; \oplus   
                                 \ar_{-\cobou_D}[d] \ar_{f}[dr]
                             & C^n  
                                 \ar^{\cobou_C}[d]\\
                               \cone f ^{n+1} & = D^{n+2}\; \oplus
                             & C^{n+1}.
               }
             \]
          \end{itemize}
          In the first case, we equip the mapping cone with the direct
          sum of the norms, in the second case, we use the maximum
          norm.
        \end{defi}

\subsubsection{Proof of the translation principle}

To prove the translation principle we just need to assemble the pieces
collected in the previous paragraphs in the right way:

\begin{bew}[Proof~(of Theorem~\ref{linkthm}).]
  The first part follows by fusing properties of mapping cones with
  the duality principle: The induced homomorphism~$H_*(f)$ is an
  isomorphism if and only if $H_*(\cone f) =0$. In view of the duality
  principle and the compatibility of mapping cones with taking the
  topological dual, this is equivalent to
  \[ 0 =     H^*\bigl( \cone f ' \bigr)
       \cong H^*\bigl( \Sigma\cone{-f'} \bigr)
       =     H^{*-1}\bigl(   \cone{-f'} \bigr); 
  \]
  the duality principle is applicable because the cone of a morphism
  of Banach chain complexes is a Banach chain complex. On the other
  hand, the~$H^{*-1}(\cone{-f'})$ are all zero if
  and only if~$H^*(-f') \colon H^*(D') \longrightarrow H^*(C')$ is an
  isomorphism. Moreover, $H^*(f') = - H^*(-f')$, and therefore the
  first part is shown. 

  For the second part, it remains to prove that $H_*(f)$ is isometric
  whenever $H^*(f')$ is an isometric isomorphism. Let $n \in \N$ and
  let $\alpha \in H_n(C)$. Using the duality principle for semi-norms
  twice, we obtain
  \begin{align*}
        \bigl\| H_n(f)(\alpha)\bigr\|
    & = \sup 
        \Bigl\{ \frac1{\supn\psi} 
        \Bigm|  \text{$\psi \in H^n(D')$ and %
                      $\krp[big]\psi{H_n(f)(\alpha)}=1$}
        \Bigr\} \\
    & = \sup 
        \Bigl\{ \frac1{\supn\psi} 
        \Bigm|  \text{$\psi \in H^n(D')$ and %
                      $\krp[big]{H^n(f')(\psi)}{\alpha}=1$}
        \Bigr\} \\
    & = \sup 
        \Bigl\{ \frac1{\supn{H^n(f')(\psi)}} 
        \Bigm|  \text{$\psi \in H^n(D')$ and %
                      $\krp[big]{H^n(f')(\psi)}{\alpha}=1$}
        \Bigr\} \\
    & = \sup
        \Bigl\{ \frac 1{\supn\varphi}
        \Bigm|  \text{$\varphi \in H^n(C')$ and %
                      $\krp\varphi\alpha = 1$}
        \Bigr\}\\
    & = \|\alpha\|. \qedhere
  \end{align*}
\end{bew}

\begin{bem}
  \emph{The converse of the second part of the translation principle
  (Theorem~\ref{linkthm}) does not hold in general:}

  Let $C=D$ be a Banach chain complex concentrated in
  degrees~$0$ and~$1$ that consists of a bounded
  operator~$\bou \colon C_1 \longrightarrow C_0$ that is not
  surjective but has dense image (e.g., the
  inclusion~$\ell^1 \hookrightarrow c_0$). In particular, the
  semi-norm on~$H_*(C) = H_*(D)$ is zero. 
  The morphism~$f \colon C
  \longrightarrow D$ given by multiplication by a constant~$c \in
  \R\setminus \{-1,0,1\}$ induces an isometric isomorphism~$H_*(f) \colon
  H_*(C) \longrightarrow H_*(D)$.

  On the other hand, the coboundary operator~$\bou' \colon C_0'
  \longrightarrow C_1'$ does not have dense image~\cite[Corollary of
  Theorem~4.12]{rudin}. Therefore, there are elements in~$H^1(D')$ of
  non-zero semi-norm. So $H^*(f')$, which is multiplication
  by~$c$, is not isometric. 
\end{bem}

\section{Isomorphisms in 
  \texorpdfstring{$\ell^1$}{l\textonesuperior}-homology}\label{isomsec} 

\noindent
In this section, we apply the translation mechanism established in the
previous section to $\ell^1$-homology, thereby gaining a uniform,
lightweight approach to proving that $\ell^1$-homology depends only on
the fundamental group (Section~\ref{appssubsec}), that
$\ell^1$-homology cannot see amenable, normal subgroups
(Section~\ref{appssubsec} and~\ref{groupsappssubsec}) and that
$\ell^1$-homology of spaces can be computed in terms of certain
projective resolutions (Section~\ref{twistedappssubsec}).

\subsection{Isomorphisms in
  \texorpdfstring{$\ell^1$}{l\textonesuperior}-homology of spaces} 
\label{appssubsec}

\noindent
We start with the simplest applications of this type, concerning
$\ell^1$-homology of spaces with $\R$-coefficients:

\begin{koro}\label{relationkoro}\nichi{02.10.2006, 03.10.2006}
  Let $f \colon (X,A) \longrightarrow (Y,B)$ be a continuous map of
  pairs of topological spaces.
  \begin{enumerate}
    \item The induced homomorphism~$\lonehom * f \colon \lonehom *
          {X,A} \longrightarrow \lonehom * {Y,B}$ is an isomorphism if
          and only if $\bch * f \colon \bch * {Y,B} \longrightarrow
          \bch * {X,A}$ is an isomorphism.
    \item If $\bch * f \colon \bch * {Y,B} \longrightarrow \bch *
          {X,A}$ is an isometric isomorphism, then $\lonehom * f$ is
          also an isometric isomorphism.
    \item In particular, $\lonehom * {X,A}$ vanishes if and only if
          $\bch * {X,A}$ vanishes.
  \end{enumerate}
\end{koro}
\begin{bew}
  By definition, $\bchn * {X,A} = (\lonech * {X,A})'$ and 
  $\bchn * {Y,B} = (\lonech * {Y,B})'$. The cochain map $\bchn * f
  \colon \bchn * {Y,B} \longrightarrow \bchn * {X,A}$ coincides
  with~$(\lonech * f)'$. Applying the translation principle
  Theorem~\ref{linkthm} to~$\lonech * f$ proves the Corollary.
\end{bew}

A discrete group~$A$ is \emph{amenable} if there is a
left-invariant mean on the set~$B(A,\R)$ of bounded functions from~$A$
to~$\R$, i.e., if there is a linear map~$m \colon B(A,\R)
\longrightarrow \R$ satisfying
\[ \fa{f \in B(A,\R)} \fa{a \in A}
   m(f) = m\bigl( b \mapsto f(a^{-1} \cdot b)
           \bigr)
\]
and
\[ \fa{f \in B(A,\R)}
        \inf\bigl\{ f(a) \bigm| a\in A \bigr\}
   \leq m(f)
   \leq \sup\bigl\{ f(a) \bigm| a\in A \bigr\}.
\]
For instance, all finite and all Abelian groups are amenable. Moreover,
the class of amenable groups is closed under taking subgroups and
quotients. An example of a non-amenable group is the free group~$\Z *
\Z$. A detailed discussion of amenability can be found in Paterson's
book~\cite{paterson}. 

\begin{koro}[Mapping theorem for $\ell^1$-homology]\label{amkerkoro}
  Let $f \colon X \longrightarrow Y$ be a continuous map between
  connected, countable \cw-complexes such that $\pi_1(f) \colon \pi_1(X)
  \longrightarrow \pi_1(Y)$ is surjective and has amenable
  kernel. Then the induced homomorphism
  \[ \lonehom * f \colon \lonehom * X \longrightarrow \lonehom * Y
  \]
  is an isometric isomorphism.
\end{koro}
\begin{bew}
  It is a classical result in the theory of bounded cohomology that in
  this situation $\bch * f \colon \bch * Y \longrightarrow \bch * X$
  is an isometric isomorphism~\cite[p.~40, Theorem~4.3]{vbc, ivanov}.
  Therefore, Corollary~\ref{relationkoro} completes the proof.
\end{bew}

Applying the mapping theorem to the classifying map~$X \longrightarrow
B\pi_1(X)$ shows in particular that the $\ell^1$-homology of a
connected, countable \cw-complex~$X$ depends only on the fundamental
group.

\subsection{Isomorphisms in
  \texorpdfstring{$\ell^1$}{l\textonesuperior}-homology of discrete
  groups} \label{groupsappssubsec}

For $\ell^1$-homology of discrete groups the translation principle
takes the following form:

\begin{koro}\label{l1groupcoro}
  Let $\varphi \colon G \longrightarrow H$ be a homomorphism of
  discrete groups, let $V$ be a Banach $G$-module, let $W$ be a Banach
  $H$-module and suppose that $f \colon V \longrightarrow \varphi^* W$
  is a morphism of Banach $G$-modules.
  \begin{enumerate}
    \item Then the homomorphism~$\lonehom * {\varphi;f} \colon
      \lonehom * {G;V} \longrightarrow \lonehom * {H;W}$ is an
      isomorphism if and only if $\bch * {\varphi;f'} \colon \bch *
      {H;W'} \longrightarrow \bch * {G;V'}$ is an isomorphism.
    \item If $\bch * {\varphi; f'}$ is an isometric isomorphism, then
      so is~$\lonehom * {\varphi;f}$.
    \item In particular, $\lonehom * {G;V} \cong \lonehom * {1;V_G}$ if
      and only if~$\bch * {G;V'} \cong \bch[big] * {1;(V')^G}$.
  \end{enumerate}
\end{koro}
\begin{bew}
  By definition, we have 
  \begin{align*}
        \lonehom * {\varphi;f}
    & = H_* \bigl( p \circ \lonech * {\varphi;f}_G\bigr), \\
        \bch * {\varphi; f'}
    & = H^* \bigl( \bchn * {\varphi;f'}^G \circ i \bigr),
  \end{align*}
  where $p \colon (\varphi^* \lonech * {H;W})_G \longrightarrow
  \lonech * {H;W} _H$ and $i \colon \bchn
  * {H;W'} ^H \longrightarrow (\varphi^* \bchn * {H;W'})^G$ denote the
  canonical projection and the inclusion respectively.

  \begin{figure}
    \begin{center}
      \begin{align*}
        \xymatrix@C=1em{%
          \bigl( \lonech * {H;W}_H \bigr)'
          \ar@{=}[r]
          \ar[d]^{p'}
          \ar@/_4pc/[dd]_{(p \circ \lonech[small] * {\varphi;f}_G)'}
        & \bigl( \lonech * {H;W}'  \bigr)^H
          \ar@{=}[r]
        & \bchn * {H;W'}^H
          \ar[d]_{i}
          \\
          \bigl( (\varphi^* \lonech * {H;W})_G \bigr)'
          \ar@{=}[r]
          \ar[d]^{(\lonech[small] * {\varphi;f}_G)'}
        & \bigl( \varphi^* \lonech * {H;W}'    \bigr)^G
          \ar@{=}[r]
        & \bigl( \varphi^* \bchn * {H;W'}      \bigr)^G
          \ar[d]_{\bchn[small] * {\varphi;f'}^G}
          \\
          \bigl( \lonech * {G;V}_G \bigr)'
          \ar@{=}[r]
        & \bigl( \lonech * {G;V}'  \bigr)^G
          \ar@{=}[r]
        & \bchn * {G;V'} ^G     
        }
      \end{align*}
    \end{center}
    \caption{Linking $\ell^1$-homology and bounded cohomology of discrete
      groups (proof of Corollary~\ref{l1groupcoro})}\label{groupstranslfig}
  \end{figure}
  A straightforward calculation shows that the diagram
  in Figure~\ref{groupstranslfig} is a commutative diagram  
  of morphisms of Banach cochain complexes, where all horizontal
  morphisms are isometric isomorphisms. 
  Thus, applying the translation principle 
  (Theorem~\ref{linkthm}) to the morphism~$p \circ \lonech *
  {\varphi;f}_G$ of Banach chain complexes proves the first 
  two parts of the corollary. 
  The third part follows because $(V_G)'$ and $(V')^G$ are
  isometrically isomorphic.
\end{bew}

An interesting consequence of the third statement is that it provides
a characterisation of amenable groups:

\begin{koro}\label{amenabilitykoro}
  For a discrete group~$G$ the following are equivalent:
  \begin{enumerate}
    \item The group~$G$ is amenable.
    \item For all Banach $G$-modules~$V$, the
      $\ell^1$-homology~$\lonehom * {G;V}$ of~$G$ with coefficients
      in~$V$ is trivial, i.e., $\lonehom * {G;V} \cong \lonehom *
      {1;V_G}$. 
  \end{enumerate}
\end{koro}
\begin{bew}
  Amenable groups can be characterised by the vanishing of bounded
  cohomology with arbitrary (dual) coefficients in non-zero
  degree~\cite{johnson, noskov}. Therefore, the claim follows with
  help of Corollary~\ref{l1groupcoro} and Section~\ref{zerosubsubsec}. 
\end{bew}

Like $\ell^1$-homology of spaces, $\ell^1$-homology of discrete groups
cannot see amenable, normal subgroups:

\begin{koro}
  Let $G$ be a discrete group, let $A \subset G$ be an amenable,
  normal subgroup, and let $V$ be a Banach $G$-module. Then the
  projection~$G \longrightarrow G/A$ induces an isometric isomorphism 
  \[ \lonehom * {G;V} \cong \lonehom * {G/A; V_A}. 
  \]
\end{koro}
\begin{bew}
 The corresponding homomorphism
  \[ \bch[norm] * {G \twoheadrightarrow G/A; V'^A \hookrightarrow V'} 
     \colon
     \bch[norm] * {G/A; {V'}^A} \longrightarrow \bch * {G;V'}
  \] 
  is an isometric isomorphism~\cite[Theorem~1,
  Corollary~8.5.2]{noskov, monod} (the case with \mbox{$\R$-co}\-ef\-fi\-cients was
  already treated by Ivanov~\cite[Section~3.8]{ivanov}). Because the
  inclusion~$V'^A \hookrightarrow V'$ is the dual of the projection~$V
  \longrightarrow V_A$, we can apply Corollary~\ref{l1groupcoro}.
\end{bew}

\begin{caveat}\label{amenablecaveat}
  Let $G$ be a discrete group and let $A \subset G$ be an amenable, 
  normal subgroup. 
  Ivanov proved that the cochain complex~$\bchn * {G/A}$ is a strong
  relatively injective $G$-resolution of the trivial
  $G$-module~$\R$~\cite[Theorem~3.8.4]{ivanov} by showing that the
  \mbox{$G$-mor}\-phisms~$\bchn * {G/A} \longrightarrow \bchn * G$
  induced by the projection~$G \longrightarrow G/A$ are split
  injective~\cite[Lemma~3.8.1 and Corollary~3.8.2]{ivanov}.

  Analogously, Park claimed that the $G$-morphisms~$\lonech n G
  \longrightarrow \lonech n {G/A}$ are split
  surjective~\cite[Lemma~2.4 and Lemma~2.5]{park} and concluded that
  the $\lonech n {G/A}$ are relatively projective
  $G$-modules. Unfortunately, Park's proof~\cite[proof of
  Lemma~2.4]{park} contains an error: the $A$-invariant mean
  on~$B(A,\R)$ provided by amenability of~$A$ in general is \emph{not}
  $\sigma$-additive.

  In fact, $\lonech n {G/A}$ in general is \emph{not} a relatively
  projective $G$-module as the following example shows: Let $G$ be an
  infinite amenable group (e.g., $G = \Z$) and $A := G$. Then the
  $G$-action on~$G/A = 1$ is trivial. However, since $G$ is infinite,
  the $G$-modules~$\lonech n G$ do not contain any non-zero
  $G$-invariant elements. Therefore, any $G$-morphism of type~$\lonech
  n {G/A} \longrightarrow \lonech n G$ must be trivial. We now
  consider the mapping problem
  \begin{align*}
    \xymatrix{%
    & \txt{\makebox[0pt]{$\lonech n {G/A} = \R$}\phantom{(\R)}} 
      \ar[d]^{\id} \ar@{..>}[dl]|{\txt{?}} 
    & \\
      \lonech n G \ar[r]_-{\pi}
    & \R \ar[r]
    & 0
    }
  \end{align*}
  with the $G$-morphism~$\pi$ given by $g_0 \cdot
  [g_1|\dots|g_n] \longmapsto 1$, which obviously admits a
  (non-equivariant) split of norm~$1$; i.e., the morphism~$\pi$ is
  relatively projective. The argument above shows that
  this mapping problem cannot have a solution, and hence that $\lonech
  n {G/A}$ cannot be a relatively projective $G$-module.

  This problem also affects several other results of Park, e.g., her
  proof of the fact that $\ell^1$-homology depends only on the
  fundamental group~\cite[Theorem~4.1]{park} and of the equivalence
  theorem~\cite[Theorem~3.7 and~4.4]{park}.
\end{caveat}

\subsection{ \texorpdfstring{$\ell^1$}{l\textonesuperior}-Homology via 
  projective resolutions}\label{twistedappssubsec}

Ivanov proved that bounded cohomology of a topological space with
$\R$-coefficients can be computed in terms of strong relatively
injective resolutions~\cite{ivanov}. The translation principle allows
us to deduce that $\ell^1$-homology of spaces also admits such a
description in terms of homological algebra:

\begin{koro}\label{l1coeffcoro}\nichi{25.10.2006, 29.11.2006}
  Let $X$ be a countable, connected \cw-complex with fundamental
  group~$G$ and let $V$ be a Banach $G$-module. 
  \begin{enumerate}
    \item There is a canonical isometric isomorphism
          \[ \lonehom * {X;V} \cong \lonehom*{G;V}. \]
    \item If $C$ is a strong relatively projective
          resolution of~$V$, then there is a canonical isomorphism
          (degreewise isomorphism of semi-normed vector spaces)
          \[ \lonehom * {X;V} \cong H_*(C_G). 
          \]
    \item If $C$ is a strong relatively projective resolution of the
          trivial Banach $G$-module~$\R$, then there is a canonical
          isomorphism (degreewise isomorphism of semi-normed vector
          spaces)
          \[ \lonehom * {X;V} \cong H_* \bigl( (C \projten V)_G
                                        \bigr).
          \]
  \end{enumerate} 
\end{koro}

Therefore, the results of Section~\ref{groupsappssubsec} are also
valid for $\ell^1$-homology with twist\-ed coefficients and hence
provide generalisations of the results presented in
Section~\ref{appssubsec}.

\begin{caveat}\label{strongcaveat}
  Ivanov proved the corresponding theorem for bounded cohomology with
  $\R$-coefficients by verifying that~$\bchn[norm] * {\ucov X}$ is a
  strong relatively injective resolution of the trivial Banach
  $G$-module~$\R$~\cite[Theorem~2.4]{ivanov}.

  The proof that the resolution~$\bchn[norm]*{\ucov X}$ is strong
  relies heavily on the fact that certain chain maps are split
  injective. However, for the same reasons as explained in
  Caveat~\ref{amenablecaveat}, it is not possible to translate these
  arguments into the language of \mbox{$\ell^1$-chain} complexes. Hence, it
  seems impossible to prove that the chain complex~$\lonech * {\ucov
    X}$ is a \emph{strong} resolution. In particular, Park's
  proof~\cite[proof of Theorem~4.1]{park} of
  Corollary~\ref{l1coeffcoro} (with $\R$-coefficients) is not
  complete.
\end{caveat}

\begin{bew}[Proof~(of Corollary~\ref{l1coeffcoro}]
   \emph{Ad~1.}\ In order to prove the first part of
  Corollary~\ref{l1coeffcoro}, we proceed as follows:
  \begin{enumerate}
    \item We establish a connection between~$\lonech * {\ucov X;V}$ and
          the strong relatively projective resolution~$\lonech * {G;V}$.
    \item The dual of this morphism, when restricted to the
          invariants, induces an isometric isomorphism on the level of
          cohomology of the invariants~\cite[Appendix~B]{loehphd};
          this is a straightforward generalisation of Ivanov's result
          that bounded cohomology with $\R$-coefficients can be
          computed in terms of strong relatively injective
          resolutions.
    \item Finally, we apply the translation principle
          (Theorem~\ref{linkthm}) to transfer this isometric
          isomorphism back to $\ell^1$-homology.
  \end{enumerate}

  \emph{First step.}
  Park~\cite[proof of Theorem~4.1]{park} constructed the following map
  (``pre-dually'' to Ivanov's construction~\cite[proof of
  Theorem~4.1]{ivanov}):

  Let $F \subset \ucov X$ be a (set-theoretic) fundamental domain of
  the $G$-action on~$\ucov X$. In the following, the vertices
  of the standard $n$-simplex~$\ssim n$ are denoted by~$v_0, \dots,
  v_n$. For a singular simplex~$\sigma \in \map[norm]{\ssim n}{\ucov
  X}$ let $g_0(\sigma), \dots, g_n(\sigma) \in G$ be the group
  elements defined inductively by the requirement that 
  \[         g_j(\sigma)^{-1} 
       \cdot \dots
       \cdot g_1(\sigma)^{-1}
       \cdot g_0(\sigma)^{-1}
       \cdot \sigma(v_j) \in F
  \]
  for all~$j \in \{0, \dots, n\}$.
  Then the map~$\eta\colon\lonech * {\ucov X} \longrightarrow \lonech
  * G$ given by
  \begin{align*}
    \lonech n {\ucov X} & \longrightarrow \lonech n G \\
    \sigma & \longmapsto 
    g_0(\sigma) 
    \cdot \bigl[g_1(\sigma) \big| \ldots \big|\, g_n(\sigma) 
          \bigr], 
  \end{align*}
  and hence also 
  $\eta_V := \eta\projten\id_V \colon \lonech * {\ucov X;V}
     \longrightarrow \lonech * {G;V}$,  
  is a morphism of Banach $G$-chain complexes. 
  Let $(\eta_V)_G \colon \lonech *{\ucov X;V}_G \longrightarrow
  \lonech * {G;V}_G$ denote the morphism of Banach chain complexes
  induced by~$\eta_V$. 

  We verify now that a different choice of fundamental domain~$F^*
  \subset \ucov X$ leads to a map chain homotopic to~$(\eta_V)_G$: By
  the fundamental lemma of homological algebra in the context of
  Banach $G$-modules~\cite[Proposition~A.7]{loehphd}, there is up to
  $G$-homotopy only one \mbox{$G$-mor}\-phism~$\lonech * {\ucov X}
  \longrightarrow \lonech * {G}$; in fact, $\lonech * {\ucov X}$ is a
  Banach \mbox{$G$-chain} complex consisting of relatively projective
  $G$-modules~\cite[p.~611]{park} and $\lonech * G$ is a strong
  relatively projective resolution
  of~$\R$~\cite[Proposition~2.19]{loehphd}. But $\eta$ and $\eta^*$,
  the map obtained via~$F^*$, are such $G$-morphisms and hence are
  $G$-homotopic. Therefore, also $\eta \projten \id_V$ and $\eta^*_V
  := \eta^* \projten \id_V$ must be $G$-homotopic, which implies that
  the induced maps~$(\eta_V)_G$ and $(\eta_V^*)_G$ are homotopic. In
  particular,
  $H_*\bigl((\eta_V)_G\bigr)
     \colon H_*\bigl( \lonech * {\ucov X;V}_G \bigr)
     \longrightarrow 
     H_*\bigl( \lonech * {G;V}_G \bigr)
  $
  does not depend on the choice of fundamental domain.

  \emph{Second step.}
  The dual of the $G$-morphism~$\eta_V$ coincides under the natural
  isometric isomorphisms~$(\lonech * {\ucov X;V})' \cong \bchn[norm] *
  {\ucov X;V'}$ and~$(\lonech * {G;V})' \cong \bchn * {G;V'}$ of
  Banach \mbox{$G$-co}chain complexes with~$\vartheta_{V'} \colon \bchn * {G;V'}
  \longrightarrow \bchn[norm] * {\ucov X;V'}$, the morphism of Banach
  $G$-cochain complexes given by
  \begin{align}\label{thetaeq}
  \begin{split}
     \bchn n {G;V'}
   & \longrightarrow \bchn[norm] n {\ucov X;V'}\\
     f 
   & \longmapsto
     \bigl( \sigma \mapsto
            f(g_0(\sigma), \dots, g_n(\sigma))
     \bigr).
  \end{split}
  \end{align}
  In other words, the diagram in Figure~\ref{twistedfig}(a) 
  is commutative. Taking $G$-invariants of this diagram yields the
  commutative diagram of morphisms of Banach cochain
  complexes depicted in Figure~\ref{twistedfig}(b).
  \begin{figure}
    \begin{tabularx}{\linewidth}{XX}
    \begin{center}
        $\xymatrix@=2em{%
          \bigl(\lonech * {G;V}\bigr)'
          \ar@{=}[d]
          \ar[r]^-{(\eta_V)'}
        & \bigl(\lonech * {\ucov X;V}\bigr)'
          \ar@{=}[d]
        \\
          \bchn * {G;V'}
          \ar[r]_-{\vartheta_{V'}}
        & \bchn[norm] * {\ucov X;V'}}
        $
    \end{center}
    &
    \begin{center}
        $\xymatrix@=2em{%
          \bigl(\lonech * {G;V} _G \bigr)'
          \ar@{=}[d]
          \ar[r]^-{(\eta_V)_G{}'}
        & \bigl(\lonech * {\ucov X;V}_G\bigr)'
          \ar@{=}[d]
          \\
          \bigl(\lonech * {G;V}'\bigr)^G
          \ar@{=}[d]
          \ar[r]^-{(\eta_V)'{}^G}
        & \bigl(\lonech * {\ucov X;V}'\bigr)^G
          \ar@{=}[d]
          \\
          \bchn * {G;V'}^G
          \ar[r]_-{(\vartheta_{V'})^G}
        & \bchn[norm] * {\ucov X;V'}^G    
        }$
    \end{center}
        \\
        \centering (a) 
      & \centering (b)
    \end{tabularx}
    \caption{Relating the morphisms~$\eta_V$ and~$\vartheta_{V'}$}\label{twistedfig}
  \end{figure}

  The restriction~$(\vartheta_{V'})^G$ to the subcomplexes of
  $G$-invariants induces an isometric isomorphism on the level of
  cohomology~\cite[Appendix~B]{loehphd}. Hence, also
  the top row of the diagram (i.e, $(\eta_V)_G{}'$) must induce an
  isometric isomorphism on the level of cohomology.

  \emph{Third step.}
  Therefore, we can derive from the translation principle
  (Theorem~\ref{linkthm}) that
  $(\eta_V)_G
     \colon
       \lonech * {X;V} 
     = \lonech * {\ucov X;V} _G 
     \longrightarrow 
     \lonech * {G;V} _G
  $ 
  induces a (canonical) isometric isomorphism on the level of
  homology. This finishes the proof of the first part. 

 \emph{Ad~2.\ and~3.}\ 
  These statements follow from the first part combined with
  the corresponding results on $\ell^1$-homology of discrete groups
  (Theorem~\ref{discreteresthm}). 
\end{bew}

For example, using this description of $\ell^1$-homology via
projective resolutions, one can construct a ``straightening'' on the
$\ell^1$-chain complex of countable, connected
\mbox{\cw-com}\-plexes~\cite[Section~4.4]{loehphd}, generalising the classical
straightening of Thurston~\cite[p.~6.3]{thurston} in the presence of
non-positive curvature. An important aspect of this generalised
straightening is that it allows to get control of the semi-norm in
measure homology~\cite[Appendix~D]{loehphd}, thereby obtaining
homological (and hence a bit more transparent) versions of the
original proofs~\cite[Section~4.3, Theorem~1.1 and~1.2]{strohm,loehmh}
that measure homology and singular homology are \emph{isometrically}
isomorphic.

\section{Simplicial volume of non-compact manifolds}\label{fincritsec}

\noindent
The definition of simplicial volume can be adapted to cover also
non-compact manifolds. In this section, we demonstrate how to utilise
$\ell^1$-homology and the results established in Section~\ref{isomsec}
to study the simplicial volume of non-compact manifolds: We first
express the simplicial volume of non-compact manifolds in terms of
$\ell^1$-homology (Section~\ref{ncsvdefsubsec}). In
Section~\ref{fincritsubsec}, we present a finiteness criterion for the
simplicial volume of non-compact manifolds. Applications of
this finiteness criterion are discussed in
Section~\ref{fincritappssubsec}.

\subsection{Simplicial volume -- the non-compact case}\label{ncsvdefsubsec}

The $\ell^1$-norm on the singular chain complex admits an obvious
extension to the chain complex of locally finite chains (notice
however, that there are locally finite chains with \emph{infinite}
$\ell^1$-norm). In particular, there is also a notion of simplicial
volume for non-compact manifolds:

\begin{defi}[Simplicial volume of non-compact manifolds]
  Let $M$ be an oriented, connected (possibly non-compact) $n$-manifold
  without boundary. Then the \emph{simplicial volume of~$M$} is
  defined by
  \[    \sv M
     := \inf\,\bigl\{ \lone c
              \bigm|  \text{$c \in \lfsing n M$ locally finite 
                      $\R$-fundamental cycle of~$M$}
              \bigr\}
     \in [0,\infty].
    \qedhere
  \]
\end{defi}

By definition, the simplicial volume of non-compact manifolds is
invariant under proper homotopy equivalences. We now provide a
description of the simplicial volume for not necessarily compact
manifolds in terms of $\ell^1$-homology:

\begin{defi}
  If $M$ is an oriented, connected $n$-manifold without boundary, we
  write~$\lonefcl M \subset \lonehom n M$ for the set of all
  homology classes in~$\lonehom n M$ that are represented by at least
  one locally finite fundamental cycle (with finite $\ell^1$-norm).
\end{defi}

If $M$ is compact, then the set~$\lonefcl M$ contains exactly one
element, namely the class~$H_n(\csing n M \hookrightarrow \lonech n
M)(\fcl M)$. However, if $M$ is non-compact, the set~$\lonefcl M$
may be empty (this happens if and only if $\sv M = \infty$) or consist
of more than one element.

\begin{prop}\label{ncsvprop}
  Let $M$ be an oriented, connected $n$-manifold without boundary.
  \begin{enumerate}
    \item Then 
          $ \sv M = \inf \bigl\{ \lone \alpha 
                          \bigm| \alpha \in \lonefcl M
                                 \subset \lonehom n M 
                          \bigr\}.
          $
    \item If $\lonehom n M = 0$, then $\sv M \in \{0,\infty\}$. 
  \end{enumerate}
\end{prop}

\begin{bew}
  The second part is an immediate corollary of the first one. We now
  prove the first part:  
  Let $j \colon \lfsing * M \cap \lonech * M \hookrightarrow
  \lonech * M$ denote the inclusion. By definition, 
  \begin{align*} 
    \sv M   = \inf\bigl\{ \lone \alpha
                  \bigm|  \alpha \in H_*(j)^{-1}(\lonefcl M)
                  \bigr\}.
  \end{align*}
  The sequence
  $                   \csing * M 
      \hookrightarrow \lfsing * M \cap \lonech * M
      \hookrightarrow \lonech * M
  $
  of inclusions of normed chain complexes shows that the middle
  complex is a dense subcomplex of the $\ell^1$-chain complex~$\lonech
  * M$. Thus, the induced map~$H_*(j) \colon
      H_*\bigl(\lfsing * M \cap \lonech * M\bigr)
      \longrightarrow
      \lonehom * M
  $
  on homology is isometric (Proposition~\ref{denseprop}). This yields 
  the desired description of~$\sv M$.
\end{bew}

For example, if $M$ is an oriented, connected manifold (of non-zero
dimension) without boundary and amenable fundamental group, then~$\sv
M \in \{0,\infty\}$.

Using the duality principle for semi-norms one also obtains a
corresponding result expressing the simplicial volume of non-compact
manifolds via bounded cohomology; however, this description is not
as convenient as the one in terms of $\ell^1$-homology.

\subsection{A finiteness criterion}\label{fincritsubsec}

In general, the simplicial volume of non-compact manifolds is not
finite -- it can even then be infinite if the manifold in question is
the interior of a compact manifold with boundary. In this case,
$\ell^1$-homology gives a necessary and sufficient finiteness
condition:

\begin{satz}[Finiteness criterion]\label{fincritthm}
  Let $(W,\bou W)$ be an oriented, compact $n$-mani\-fold
  with boundary~$\bou W$. Then the following
  are equivalent:
  \begin{enumerate}
    \item The simplicial volume of the interior~$W^\circ$ is finite. 
    \item The manifold~$\bou W$ is \emph{$\ell^1$-invisible,} i.e.,
      \[ H_{n-1}\bigl( \csing * {\bou W} \hookrightarrow \lonech * {\bou W}
                \bigr) 
                     (\fcl{\bou W}) = 0
         \in \lonehom{n-1}{\bou W}. 
      \]
  \end{enumerate}
\end{satz}

In particular, by combining this finiteness criterion with
Proposition~\ref{denseprop}, we obtain Gromov's necessary
condition~\cite[p.~17]{vbc}: If $\sv{W^\circ} < \infty$, then $\sv
{\bou W} = 0$. Notice that in contrast to Gromov's estimate
of the simplicial volume by the minimal volume~\cite[p.~12,
p.~73]{vbc}, the finiteness criterion is purely topological and can 
be proved by elementary means.

While it is clear that every $\ell^1$-invisible manifold has vanishing
simplicial volume by Proposition~\ref{denseprop}, it is an open
problem whether every oriented, closed, connected manifold with
vanishing simplicial volume is already $\ell^1$-invisible. 

Because the evaluation map linking bounded cohomology and
$\ell^1$-homology is continuous, bounded cohomology can detect only
whether the semi-norm of a given class in $\ell^1$-homology is zero,
but not if the class itself is zero. Therefore, the finiteness
criterion as stated above cannot be phrased in terms of bounded
cohomology.

\begin{bew}[Proof~(of Theorem~\ref{fincritthm}).]
  The theorem trivially holds if the boundary~$\bou W$ is empty;
  therefore, we assume for the rest of the proof that~$\bou W \neq
  \emptyset$. The homeomorphism~\cite{brownflat, connelly}
  \[ W^\circ \cong W \sqcup_{\bou W} \bou W \times [0, \infty) 
             =:    M
  \]
  shows that we can look at the notationally more convenient
  manifold~$M$ instead of~$W^\circ$.

   \begin{list}{}
    {\renewcommand{\makelabel}[1]{\mbox{#1}\hfil}%
     \settowidth{\labelwidth}{$2 \Rightarrow 2$}%
     \setlength{\leftmargin}{\labelwidth+\labelsep}%
     \setlength{\listparindent}{\parindent}%
     \setlength{\itemsep}{0pt}%
     \setlength{\parsep}{0pt}%
    }
    \item[$1 \Rightarrow 2$] 
      Suppose that the simplicial volume~$\sv{W^\circ} = \sv M$ is
      finite.  In other words, there is a locally finite fundamental
      cycle~$c = \sum_{j\in\N} a_j \cdot \sigma_j \in \lfsing n {M}$ of~$M$ with $\lone{c} < \infty.$ 
      We now restrict~$c$ to a cylinder lying in~$\bou W \times [0,
      \infty) \subset M$. The boundary of this restriction is a
      fundamental cycle of~$\bou W$ and the restriction itself gives
      rise to the desired boundary in the $\ell^1$-chain complex: 

      More precisely, for $t \in (0, \infty)$ we consider the
      cylinder~$Z_t := \bou W \times [t, \infty)$  
      and the projections~$p_t \colon \bou W \times [0,
      \infty) \longrightarrow Z_t$ and~$q_t \colon \bou W \times [0,
      \infty) \longrightarrow \bou W \times [0,t]$;  the
      notation is illustrated in Figure~\ref{onetwofig}.
      \begin{figure}
        \begin{center}
          \includegraphics{fincritpics1.mps}
        \end{center}
        \caption{The proof of ``$1 \Rightarrow 2$'' of the finiteness
                 criterion}\label{onetwofig}
      \end{figure}

      Because $c$ is locally finite, there exists a~$t \in (0,\infty)$
      such that the restriction~$c|_{Z_t} \in \lfsing n M$ of~$c$
      to~$Z_t$ does not meet~$W$; by definition, $c|_{Z_t} =
      \sum_{j\in J_t} a_j \cdot \sigma_j$, where $J_t := \{ j \in \N
      \mid \im \sigma_j \cap Z_t \neq \emptyset\}$.  It is not
      difficult to see that the chain~$\lfsing n {p_t}(c|_{Z_t})$ is a relative
      fundamental cycle of~$(Z_t, \bou W \times \{t\})$ and hence that
      $z_t := \bou(\csing n {p_t}(c|_{Z_t}))$ is a fundamental cycle of~$\bou W
      \times \{t\}$.  On the other hand, $\lone c$ is finite, so
      \[ b_t := \lonech n  {q_t}(z_t)
             \in \lonech[big] n {\bou W \times \{t\}}.
      \]
      By construction, $\bou b_t = z_t$, which proves that $\bou W \times
      \{t\}$ is $\ell^1$-invisible. Hence, $\bou W$ is also
      $\ell^1$-invisible. 

    \item[$2 \Rightarrow 1$]
      Conversely, suppose that part~2 is satisfied, i.e., that $\bou
      W$ is $\ell^1$-invisible. Therefore, there is a~$b \in \lonech n
      {\bou W}$ such that $z:= - \bou b$ is a fundamental cycle
      of~$\bou W$. Adding the boundaries of the partial sums~$(\sum_{j
      = 0}^{k-1} b_j)_{k
        \in \N} \subset \lonech n {\bou W}$ of~$b$ to~$z$ yields 
       a sequence of
      fundamental cycles~$(z_k)_{k \in \N} \subset \csing{n-1}{\bou
        W}$ of~$\bou W$ and a sequence of chains~$(b_k)_{k \in \N}
      \subset \csing n {\bou W}$ satisfying
      \begin{align*}
        \fa{k \in \N} \bou b_k     & = z_{k+1} - z_k, \\
        \sum_{k \in \N} \lone{b_k} & < \infty. 
      \end{align*}
      Moreover, $\lim_{k \rightarrow \infty} \lone[norm]{z_k} =
      0$. Thus, by choosing a suitable subsequence of~$(z_k)_{k \in
        \N}$ we can even find two such sequences such that
      additionally
      \[ \sum_{k \in \N} \lone{z_k} < \infty
      \]
      holds~\cite[Proposition~6.4]{loehphd}. Now the idea is --
      similarly to Gromov's argument in a special
      case~\cite[p.~8]{vbc} -- to take a relative fundamental cycle
      of~$(W, \bou W)$ and to glue the~$(b_k)_{k \in \N}$ to its
      boundary. To ensure that the resulting chain is locally finite,
      we spread out the chain~$\sum_{k \in\N} b_k$ over the
      cylinder~$\bou W \times [0,\infty)$.

      More precisely, let~$c \in \csing n {W}$ be a relative
      fundamental cycle of the manifold~$(W,\bou W)$ with
      boundary. Then $\bou c \in \csing {n-1}{\bou W}$ is
      a fundamental cycle of the oriented, compact manifold~$\bou
      W$. Of course, we may assume that $\bou c = z_0$.
      
      \begin{figure}
        \begin{center}
          \includegraphics{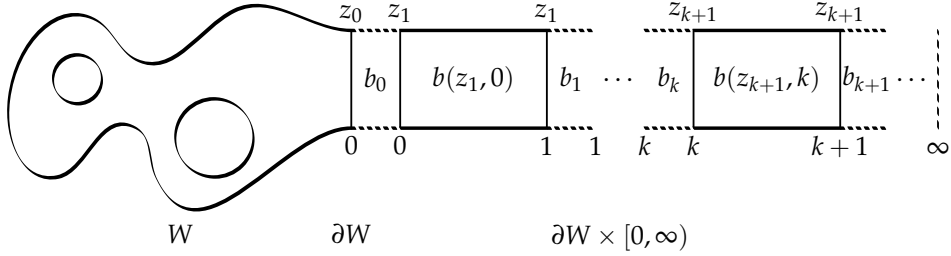}
        \end{center}
        \caption{The proof of ``$2 \Rightarrow 1$'' of the finiteness
          criterion}\label{twoonefig}
      \end{figure}
      The spreading out of~$(b_k)_{k \in \N}$ is achieved by using the
      following chains: For any cycle~$z \in \csing{n-1}{\bou W}$ and
      $k \in \N$ we can find a chain~$b(z,k) \in \csing n {\bou W
      \times [0, \infty)}$ such that
      \[        \bou\bigl(b(z,k)\bigr) 
           =    \csing{n-1}{j_{k+1}} (z)
              - \csing{n-1}{j_k}(z),  
           \quad
                \lone[big]{b(z,k)}
           \leq n \cdot \lone z ;
      \]
      here, $j_k \colon \bou W \hookrightarrow \bou W \times \{k\}
      \hookrightarrow \bou W \times [0, \infty)$ denotes the
      inclusion. For example, such a chain~$b(z,k)$ can be constructed
      by looking at the canonical triangulation of~$\ssim {n-1} \times
      [0,1]$ into $n$-simplices. 
      We set (see also Figure~\ref{twoonefig})
      \[    b 
         := \sum_{k \in \N} 
            \bigl( \csing n {j_k} (b_k) + b(z_{k+1}, k) 
            \bigr)  
      \]
      and $\overline c := c + b$. Because all~$b_k$ and all~$b(z_{k+1},
      k)$ are finite, the stretched chain~$b$ is a well-defined
      locally finite $n$-chain of~$M$. Therefore, also~$\overline c
      \in \lfsing n {M}$. By construction, $\overline c$ is a cycle
      and $\overline c|_{W^\circ} = c |_{W^\circ}$; hence, $\overline c$
      is a locally finite fundamental cycle of~$M$. Furthermore, 
      $\lone{\overline c} \leq \lone c + \lone b$,
      which shows that $\sv M < \infty$. \qedhere
   \end{list}
\end{bew}

\subsection{Applications}\label{fincritappssubsec}

Before discussing applications of the finiteness criterion
(Theorem~\ref{fincritthm}), we first have a tour through the zoo of
$\ell^1$-invisible manifolds:

\begin{bsp}[$\ell^1$-Invisibility]\label{zooex}\hfil
  \begin{itemize}
    \item \emph{Vanishing $\ell^1$-homology.} 
          By definition, any oriented, closed $n$-manifold~$M$
          satisfying $0 = \lonehom n M = \lonehom n {\pi_1(M)}$ is 
          $\ell^1$-invisible. In particular, manifolds with amenable
          fundamental group are $\ell^1$-invisible.
    \item \emph{Vanishing bounded cohomology.} Moreover, any oriented, closed
          $n$-manifold with $0 = \bch {n+1} M = \bch {n+1}{\pi_1(M)}$ is
          $\ell^1$-invisible; this follows from the fact that such
          manifolds satisfy the so-called \emph{uniform boundary
            condition in degree~$n$}~\cite[Theorem~2.8,
          Proposition~6.8]{mm, loehphd}. However, not all
          $\ell^1$-invisible $n$-manifolds satisfy the uniform boundary
          condition in degree~$n$~\cite[Example~6.9]{loehphd}.  
   \item \emph{Functoriality.} 
          Clearly, if $M \longrightarrow N$ is a continuous map of
          non-zero degree between oriented, closed manifolds of the
          same dimension and if $M$ is $\ell^1$-invisible, then so
          is~$N$. 

          Similarly, if the oriented, closed, connected
          $n$-manifold~$M$ admits a self-map~$f$ with $|\deg(f)| \geq
          2$, then $M$ is $\ell^1$-invisible: Let $z \in \csing n M$
          be a fundamental cycle of~$M$ and let $b \in \csing{n+1}M$
          with $\bou b = z - 1/{\deg f} \cdot \csing n f (z)$.  Then
          \[ \overline b := \sum_{k \in \N}
                              \frac{1}{(\deg f) ^k}
                              \cdot \csing {n+1} f ^k (b) 
          \]
          lies in~$\lonech {n+1}M$ and $z = \bou \overline b$, i.e.,
          $M$ is $\ell^1$-invisible. 
   \item  \emph{Products.} 
          If $M$ and~$N$ are oriented, closed, connected manifolds,
          and if $M$ is \mbox{$\ell^1$-in}\-visible, then using the
          $\ell^1$-version of the homological cross product on
          singular chains shows that also the product~$M \times N$ is
          $\ell^1$-invisible. 
   \item \emph{Gluings.} 
          Let $M$ and $N$ be oriented, closed, connected,
          $\ell^1$-invisible manifolds of 
          the same dimension at least~$3$. Then the connected sum~$M
          \consum N$ is also $\ell^1$-invisible:

          Let $j_M \colon M \longrightarrow M \lor N$ and $j_N \colon
          N \longrightarrow M \lor N$ be the inclusions. The 
          Mayer-Vietoris sequence for~$M \lor N$ shows that in
          non-zero degree $\hsing *
          {j_M} \oplus \hsing * {j_N}$ is an isomorphism
          mapping~$(\fcl M, \fcl N)$ to~$\fcl{M\consum N}$. Because
          $M$ and $N$ are $\ell^1$-invisible, the lowest horizontal
          map in Figure~\ref{zoofig}(a) maps~$(\fcl M, \fcl N)$
          to~$0$. 

          On the other hand, the pinching map~$f \colon M \consum N
          \longrightarrow M \lor N$ induces an isomorphism on the
          level of fundamental groups and hence induces an isomorphism
          in~$\ell^1$-homology (Corollary~\ref{amkerkoro}). Therefore,
          we can read off the commutative diagram in
          Figure~\ref{zoofig}(a) that $M \consum N$ is $\ell^1$-invisible.

          More generally, the class of $\ell^1$-invisible manifolds of
          dimension at least~$3$ is also closed under amenable
          gluings~\cite[Proposition~6.10]{loehphd}. 
          \begin{figure}
            \begin{center}
            \begin{tabularx}{\linewidth}{@{}XX@{}}
               \begin{center}\makebox[0pt]{$
                  \xymatrix{%
                     \hsing * {M \consum N}
                     \ar[r]
                     \ar[d]_{\hsing m f}
                   & \lonehom * {M \consum N} 
                     \ar[d]^{\lonehom m f}_{\cong}
                   \\
                     \hsing * {M \lor N}
                     \ar[r]
                   & \lonehom * {M \lor N} 
                   \\
                     \hsing * M \oplus \hsing m N
                     \ar[r]
                     \ar[u]^{\hsing * {j_M} \oplus \hsing * {j_N}}
                   & \lonehom * M \oplus \lonehom * N
                     \ar[u]_{\lonehom * {j_M} \oplus \lonehom * {j_N}}
                  }
               $}\end{center}
              & \begin{center}$
                  \xymatrix{%
                      \hsing * M \ar[r] 
                      \ar[d]_{\hsing * p} 
                    & \lonehom * M
                      \ar[d]^{\lonehom * p}_{\cong}
                    \\
                      \hsing * B \ar[r]
                    & \lonehom * B
                  }
               $\end{center}
                \\
                \centering (a)
              & \centering (b)
            \end{tabularx}
            \end{center}
            \caption{Proof of Example~\ref{zooex}}\label{zoofig}
          \end{figure}

   \item  \emph{Fibrations.}
          If $p \colon M \longrightarrow B$ is a fibration of
          oriented, closed, connected manifolds whose fibre~$F$ is
          also an oriented, closed, connected manifold of non-zero
          dimension and if~$\pi_1(F)$ is amenable, then $M$ is
          $\ell^1$-invisible:

          A spectral sequence argument yields~$\dim B \leq \dim M -1$.
          In particular, $\hsing * p (\fcl M) = 0 \in \hsing * B$.
          The long exact sequence of homotopy groups associated with
          the fibration~$p$ shows that $\pi_1(p)$ is surjective and
          that the kernel of~$\pi_1(p)$ is a homomorphic image of the
          amenable group~$\pi_1(F)$; thus, $\ker \pi_1(p)$ is
          amenable~\cite[Proposition~1.12 and~1.13]{paterson}.
          Therefore, $\lonehom * p$ is an isometric isomorphism
          (Corollary~\ref{amkerkoro}), and we deduce from
          Figure~\ref{zoofig}(b) that $M$ is $\ell^1$-invisible.
   \item \emph{Circle actions.} If $M$ is a smooth, oriented, closed
          manifold admitting a smooth $S^1$-action that is either free
          or has at least one fixed point, then $M$ is
          \mbox{$\ell^1$-in}\-visible: 

          In the first case, we can apply the same argument as for
          fibrations with amenable fibres because~$\pi_1(S^1) \cong
          \Z$ is amenable. 
          
          In the second case, it is known that the map on singular
          homology induced by the classifying map~$M \longrightarrow
          B\pi_1(M)$ maps~$\fcl M$ to~$0$~\cite[p.~95,
          Lemma~1.42]{vbc,lueck}, and by Corollary~$\ref{amkerkoro}$,
          the classifying map induces an isometric isomorphism
          on~$\ell^1$-homology.
    \item \emph{Proportionality.}
         If $M$ and $N$ are smooth, oriented, closed, connected manifolds
         equipped Riemannian metrics such that the Riemannian universal
         coverings of~$M$ and~$N$ are isometric, then $M$ is
         $\ell^1$-invisible if and only if $N$ is
         \mbox{$\ell^1$-in}\-visible~\cite[Proposition~6.10]{loehphd}; the proof
         of this fact is based on an~$\ell^1$-version of measure
         homology. 
    \item \emph{Relation with curvature.}
          Let $M$ be an oriented, closed, connected Riemannian manifold. 
          \begin{itemize}
            \item If $M$ has positive sectional curvature, then $\pi_1(M)$ is
              finite~\cite[Theorem~11.8]{lee}, hence amenable. In particular, $M$ is
              $\ell^1$-invisible.
            \item If $M$ is flat, then $M$ is $\ell^1$-invisible by
              proportionality, because any oriented, closed, connected
              flat manifold has the same Riemannian universal covering
              as the torus of the same dimension.
            \item If $M$ has negative sectional curvature, then $\sv M
              \neq 0$~\cite{inoue} and so $M$ is not $\ell^1$-invisible.
            \qedhere
          \end{itemize}
  \end{itemize}
\end{bsp}

Equipped with this list of examples of $\ell^1$-invisible manifolds,
we apply the finiteness criterion and the description of the
simplicial volume of non-compact manifolds in terms of
$\ell^1$-homology (Proposition~\ref{ncsvprop}) to exhibit a number of
simple examples illustrating the simplicial volume of non-compact
manifolds:

  \subsubsection{Vanishing results}

    If $(W, \bou W)$ is an oriented, compact connected $n$-manifold with
    $\ell^1$-invisible boundary and $\lonehom n W =0$, then $\sv
    {W^\circ} =0$; this follows from the finiteness criterion
    (Theorem~\ref{fincritthm}) and Proposition~\ref{ncsvprop}.

    For instance, it follows that $\sv {\R^n} =0 $ for all~$n \in
    \N_{>1}$ because the sphere~$S^{n-1}$ is
    \mbox{$\ell^1$-in}\-visible.  On the other hand, the finiteness
    criterion and $\sv[norm]{S^0} = 2$ imply that $\sv\R = \infty$. In
    particular, the simplicial volume of non-compact manifolds is in
    general not invariant under homotopy equivalences that are not
    proper.

    Notice that $\sv{\Hyp^n} = \sv{\R^n} = 0$ for~$n \in \N_{>1}$
    despite of~$\Hyp^n$ being hyperbolic. On the other hand, for
    certain classes of non-compact, negatively curved manifolds of
    finite volume non-vanishing results can be proved by more advanced
    means~\cite{hilbert}.

  \subsubsection{Non-compact manifolds with finite, non-zero simplicial
      volume}

    If $M$ is an oriented, closed, connected manifold with~$\sv M \neq
    0$ of dimension~$n\geq 2$ (for example, a closed hyperbolic
    $n$-manifold), and if $N$ is a non-compact manifold obtained
    from~$M$ by removing a finite number of points, then
    $0 < \sv N < \infty$.
    
    This can be seen as follows: By construction, $N$ is the interior
    of a compact manifold~$(N', \bou N')$ whose boundary is a disjoint
    union of~$(n-1)$-spheres. Because $S^{n-1}$ is $\ell^1$-invisible,
    the finiteness criterion (Theorem~\ref{fincritthm}) yields~$\sv N
    < \infty$.

    Why is $\sv N$ non-zero? A straightforward computation shows that
    $\sv N \geq \sv{N',\bou N'}$
    holds~\cite[Proposition~5.12]{loehphd}, where $\sv{N',\bou N'}$ is
    the infimum of the $\ell^1$-norms of all relative fundamental
    cycles of~$(N',\bou N')$. Using the fact that $D^n$ satisfies
    the uniform boundary condition in degree~$n -
    1$~\cite[Theorem~2.8]{mm} and that $H_{n-1}(D^n) = 0$, we find
    a~$K \in \R_{>0}$ with the following property: Every relative
    fundamental cycle~$z' \in \csing n {N'}$ of~$(N', \bou N')$ can be 
    extended to a fundamental cycle~$z \in \csing n M$ of~$M$ with
    \[ \lone{z} \leq \lone{z'} + K \cdot \lone{\bou z'}
                \leq \bigl(1 + K \cdot (n+1)\bigr) \cdot \lone{z'}.  
    \]
    Therefore, $\sv N \geq \sv{N', \bou N'} \geq 1/(1 + K \cdot(n+1))
    \cdot \sv M > 0$, as claimed.

%
%
%
    
  \subsubsection{Products}

  The simplicial volume of products of two manifolds can be estimated
  from below by the product of the simplicial volume of the factors if
  one of the factors is compact~\cite[p.~17f, Theorem~C.7]{vbc,
    loehphd}; however, in the case that the compact factor has
  vanishing simplicial volume and the other factor has infinite
  simplicial volume this estimate is inconclusive. In a special
  case, $\ell^1$-invisibility determines the outcome for such products:

    \begin{prop}
      Let $M$ be an oriented, closed, connected $n$-manifold. 
      Then
      \[ \sv{M \times \R} 
         = \begin{cases}
           0      & \text{if $M$ is $\ell^1$-invisible,}\\
           \infty & \text{otherwise.}
           \end{cases}
      \]
    \end{prop}
    \begin{bew}
      Because $M \times \R$ is homeomorphic to the interior of the
      compact manifold~$M \times [0,1]$ with boundary~$M \sqcup M$,
      the finiteness criterion (Theorem~\ref{fincritthm}) shows that
      $\sv{M \times \R}$ is finite 
      if and only if $M$ is $\ell^1$-invisible. 

      In the case that $M$ is $\ell^1$-invisible, the proof of the
      finiteness criterion provides us with a locally finite chain~$c
      \in \lfsing{n+1}{M \times [0, \infty)}$ such that the
      sequence~$(c_k)_{k \in \N}$ defined by~$c_k := c|_{M \times [k,
        \infty)} \in \lfsing{n+1}{M \times [k, \infty)}$ has the
      following properties: 
      \begin{itemize}
        \item For all~$k \in \N$ we have~$\bou c_k \in \csing n {M
            \times \{k\}}$ and $c_k$ is a relative locally finite
          fundamental cycle of the half-open cylinder~$(M \times [k,
          \infty), M \times \{k\})$. 
        \item Furthermore, $\lim_{k \rightarrow \infty} \lone{c_k} = 0$.
      \end{itemize}
      For $k \in \N$ we consider the mirror chain
      $ \overline{c_k} := \lfsing{n+1}{\id_M \times r_k}(c_k)$
      where $r_k \colon \R \longrightarrow \R$ denotes reflection
      at~$k$. Then $c_k - \overline{c_k} \in \lfsing{n+1}{M \times
      \R}$ is a locally finite fundamental cycle of~$M \times \R$ and
      $\sv{M \times \R}
         \leq \inf_{k \in \N} \,\lone{c_k - \overline{c_k}}
         \leq 2 \cdot \inf_{k \in \N} \,\lone{c_k}
         = 0.
      $
    \end{bew}

    Hence, any oriented, closed, connected manifold with vanishing
    simplicial volume that is not $\ell^1$-invisible would produce the
    first example of two manifolds $M$ and~$N$ satisfying $\sv M = 0$,
    $\sv N = \infty$ and $\sv{M\times N} \neq 0$.

    A related problem is to find an example of two non-compact
    manifolds whose product has non-zero simplicial volume. Using the 
    finiteness criterion we obtain:

    \begin{bsp}
      Let $(M, \bou M)$ be an oriented, compact, connected surface of
      genus at least~$1$ with non-empty boundary. Then $\sv{M^\circ
        \times \R} = \infty$:

      By construction, $M^\circ \times \R$
      is the interior of the compact manifold~$M \times [0,1]$ whose
      boundary is homeomorphic to~$M \consum M$ and hence 
      is an oriented, closed, connected surface of genus at least~$2$.
      Because hyperbolic manifolds are not $\ell^1$-invisible, the
      finiteness criterion shows that
      $\sv{M^\circ \times \R} = \infty$.
    \end{bsp}

\enlargethispage{-1.5cm}

\vfill

\small
\noindent
\emph{Clara L\"oh} \\[1ex]
\begin{tabular}{@{}ll}
  address:       & Fachbereich Mathematik und Informatik der
                   \wwu~M\"unster\\
                 & Einsteinstr.~62\\
                 & 48149~M\"unster, Germany \\[1ex]
  email:         & \textsf{clara.loeh@uni-muenster.de}\\
  \smaller{URL}: & \textsf{http://wwwmath.uni-muenster.de/u/clara.loeh}
\end{tabular}
\end{document}